\documentclass[aos,MSNbibl,seceqn,citesort,dvips]{arximspdf}
\usepackage{euscript}
\usepackage{dcolumn}
\usepackage{graphicx}

\doi{10.1214/12-AOS970} 
\volume{40}
\issue{2}
\pubyear{2012}
\firstpage{694}
\lastpage{726}

\makeatletter

\newcolumntype{d}[1]{D{.}{.}{#1}}

\newtheorem{theorem}{Theorem}

\newtheorem{lemma}{Lemma}

\newtheorem{corollary}{Corollary}

\newtheorem{proposition}{Proposition}

\newproclaim{Remark}{Remark}
\newproclaim{Example}{Example}

\renewcommand{\hat}{\widehat}
\newcommand{\wh}{\widehat}
\newcommand{\wt}{\widetilde}

\newcommand{\WN}{\mathrm{WN}}
\newcommand{\cov}{\operatorname{Cov}}
\newcommand{\de}{\delta}
\newcommand{\la}{\lambda}
\newcommand{\ve}{{\varepsilon}}

\newcommand{\bA}{{\mathbf A}}
\newcommand{\bB}{{\mathbf B}}
\newcommand{\bF}{{\mathbf F}}
\newcommand{\bE}{{\mathbf E}}
\newcommand{\bG}{{\mathbf G}}
\newcommand{\bH}{{\mathbf H}}
\newcommand{\bI}{{\mathbf I}}
\newcommand{\bM}{{\mathbf M}}
\newcommand{\bP}{{\mathbf P}}
\newcommand{\bS}{{\mathbf S}}
\newcommand{\bW}{{\mathbf W}}
\newcommand{\bX}{{\mathbf X}}
\newcommand{\bu}{{\mathbf u}}
\newcommand{\by}{{\mathbf y}}
\newcommand{\bSigma}{\bolds{\Sigma}}
\newcommand{\bgamma}{\bolds{\gamma}}
\newcommand{\bga}{\bolds{\gamma}}
\newcommand{\bve}{\bolds\varepsilon}
\newcommand{\bmu}{\bolds{\mu}}
\newcommand{\bC}{{\mathbf C}}
\newcommand{\bD}{{\mathbf D}}
\newcommand{\bzero}{{\mathbf0}}

\newcommand{\pcon}{\stackrel{P}{\longrightarrow}}

\newcommand{\calM}{{\mathcal M}}
\newcommand{\RR}{\EuScript R}
\newcommand{\norm}[1]{\|#1\|}

\setattribute{abstract}{width}{290pt}

\makeatother

\begin{document}
\begin{frontmatter}

\title{Factor modeling for high-dimensional time series: Inference for the number of factors\thanksref{T1}}
\runtitle{Factor modeling for high-dimensional time series}

\thankstext{T1}{Supported in part by
the Engineering and Physical Sciences Research Council
of the United Kingdom.}

\begin{aug}
\author[A]{\fnms{Clifford} \snm{Lam}\ead[label=e1]{c.lam2@lse.ac.uk}}
\and
\author[A]{\fnms{Qiwei} \snm{Yao}\corref{}\ead[label=e2]{q.yao@lse.ac.uk}}
\runauthor{C. Lam and Q. Yao}
\affiliation{London School of Economics and Political Science, and
London School of Economics and Political Science and
Guanghua School of Management, Peking~University}
\address[A]{Department of Statistics\\
London School of Economics\\
London, WC2A 2AE\\
United Kingdom\\
\printead{e1}\\
\phantom{E-mail: } \printead*{e2}} 
\end{aug}

\received{\smonth{3} \syear{2011}}
\revised{\smonth{1} \syear{2012}}

%
\begin{abstract}
This paper deals with the factor modeling for high-dimensional time series
based on a dimension-reduction viewpoint.
Under stationary settings, the inference is simple in the sense that
both the number of factors and the factor loadings are estimated in
terms of an eigenanalysis for a~nonnegative definite matrix, and is
therefore applicable when the dimension of time series is on the order
of a few thousands. Asymptotic properties of the proposed method are
investigated under two settings: (i) the sample size goes to infinity
while the dimension of time series is fixed; and (ii) both the sample
size and the dimension of time series go to infinity together. In
particular, our estimators for zero-eigenvalues enjoy faster
convergence (or slower divergence) rates, hence making the estimation
for the number of factors easier. In particular, when the sample size
and the dimension of time series go to infinity together, the
estimators for the eigenvalues are no longer consistent. However, our
estimator for the number of the factors, which is based on the ratios
of the estimated eigenvalues, still works fine. Furthermore, this
estimation shows the so-called ``blessing of dimensionality'' property
in the sense that the performance of the estimation may improve when
the dimension of time series increases. A two-step procedure is
investigated when the factors are of different degrees of strength.
Numerical illustration with both simulated and real data is also
reported.
\end{abstract}

%

\end{frontmatter}

\section{Introduction}
\label{sec1}

The analysis of multivariate time series data is of increased
interest and importance in the modern information age. Although
the\vadjust{\goodbreak}
methods and the associate theory for univariate time series analysis
are well developed and understood, the picture for the multivariate
cases is less complete. In spite of the fact that the conventional
univariate time series models (such as ARMA) and the associated
time-domain and frequency-domain methods have been formally extended
to multivariate cases, their usefulness is often limited. One may
face serious issues such as the lack of model identification {{or}}
flat likelihood functions. In fact vector ARMA models are seldom
used directly in practice. Dimension-reduction via, for example,
reduced-rank structure, structural indices, scalar component models
and canonical correlation analysis is more pertinent in modeling
multivariate time series data. See~\cite{H70,P81,L93,R97}.

In this paper we deal with the factor modeling for multivariate time
series from a dimension-reduction viewpoint.
Differently from the factor analysis for
independent observations, we search for the factors
which drive the serial dependence of the original time series.
Early attempts in this direction include
\cite{A63,PRT74,B81,SG84,PB87,SS89,PY08}.
More recent efforts focus on the inference when
the dimension of time series is as large as or even greater than the
sample size;
see, for example,~\cite{LYB11}
and the references within.
High-dimensional time series data are often encountered nowadays in many
fields including finance, economics, environmental and medical studies.
For example, understanding the dynamics of the returns of large numbers
of assets
is the key for asset pricing, portfolio allocation, and risk
management. Panel time series are commonplace in studying
economic and business phenomena. Environmental time series are often
of a high dimension due to a large number of indices monitored
across many different locations.

Our approach is from a dimension-reduction
point of view. The model adopted can be traced back at least to that of
\cite{PB87}.
We decompose a high-dimensional time series into two parts: a
dynamic part driven by, hopefully, a lower-dimensional factor time
series, and a static part which is a vector white noise. Since the
white noise exhibits no serial correlations, the decomposition is
unique in the sense that both the number of factors (i.e., the
dimension of the factor process) and the factor loading space in our
model are identifiable. Such a conceptually simple decomposition
also makes the statistical inference easy. Although the setting allows the
factor process to be nonstationary (see~\cite{PY08};
also Section~\ref{sec31} below),
we focus on stationary models only in this paper: under the stationary
condition,
the estimation for both the
number of factors and the factor loadings is carried out in an
eigenanalysis for a nonnegative definite matrix, and is therefore
applicable when the dimension of time series is on the order of a
few thousands. Furthermore, the asymptotic properties of the proposed
method are
investigated under two settings: (i) the sample size goes to
infinity while the dimension of time series is fixed; and (ii) both
the sample size and the dimension of time series go to infinity
together. In particular, our estimators for zero-eigenvalues enjoy
the faster convergence (or {{slower}} divergence) rates, from which
the proposed ratio-based estimator for the number of factors
benefits. In fact when all the factors are strong, the performance of our
estimation for the number of factors improves when
the dimension of time series increases. This phenomenon is coined as
``blessing of dimensionality.''

The new contributions of this paper include: (i) the ratio-based
estimator for
the number of factors and the associated asymptotic theory which underpins
the ``blessing of dimensionality'' phenomenon observed in numerical
experiments, and
(ii) a two-step estimation procedure when the factors are of different
degrees of
strength. We focus on the results related to the estimation for the
number of factors in this paper. The results on the estimation of the
factor loading space under the assumption that the number of factors is
known are reported in~\cite{LYB11}.

There exists a large body of literature in econometrics and finance
on factor models for high-dimensional time series.
However, most of them are based on a different viewpoint, as those
models attempt to identify the \textit{common factors} that affect the
dynamics of most original component series. In analyzing economic and financial
phenomena, it is often appealing
to separate these common factors from the so-called idiosyncratic
components: each idiosyncratic component may at most
affect the dynamics of a few original time series. An idiosyncratic
series may exhibit serial correlations and, therefore, may be a time
series itself. This poses technical difficulties in both model
identification and inference. In fact the rigorous definition of the
common factors and the idiosyncratic components can only be established
asymptotically when the dimension of time series tends to infinity;
see~\cite{CR83,FHLR00}.
Hence those
factor models are only asymptotically identifiable.
According to the definition adopted in this paper, both ``the common
factors'' and those
serially correlated idiosyncratic components will be identified as factors.
This is not ideal for the applications with the purpose to identify those
common factors. However, this makes the tasks of model identification and
inference much simpler.

The rest of the paper is organized as follows. The model and the
estimation methods are introduced in Section~\ref{sec3}. The sampling
properties of the estimation methods are investigated in Section
\ref{sec5}. Simulation results are inserted whenever appropriate to
illustrate the various asymptotic properties of the methods.
Section~\ref{sec6} deals with the cases when different factors are of
different strength, for which a two-step estimation procedure is
preferred. The analysis of two real data sets is reported in Section
\ref{sec7}. All mathematical proofs are relegated to the
\hyperref[app]{Appendix}.

\section{Models and estimation}
\label{sec3}

\subsection{Models}
\label{sec31}

If we are interested in the linear dynamic structure of $\by_t$ only,
conceptually we may think that $\by_t$ consists of two parts: a static
part (i.e., a white noise), and a\vadjust{\goodbreak} dynamic component driven by,
hopefully, a low-dimensional process.
This leads to the decomposition:
%
\begin{equation} \label{c1}
\by_t = \bA{\mathbf x}_t + \bve_t,
\end{equation}
where ${\mathbf x}_t $ is an $r\times1$ latent process with (unknown)
$r\le p$,
$\bA$ is a $p\times r$ unknown constant matrix,
and $\bve_t \sim\WN( \bmu_\ve,   \bSigma_\ve)$ is a vector
white-noise process.
When $r$ is much smaller than $p$, we achieve an effective
dimension-reduction, as then the serial
dependence of $\by_t$ is driven by that of a much lower-dimensional
process ${\mathbf x}_t$.
We call ${\mathbf x}_t$ a factor process.
The setting~(\ref{c1}) may be traced back
at least to~\cite{PB87};
see also its further development
in dealing with cointegrated factors in~\cite{PP06}.

Since none of the elements on the RHS of~(\ref{c1}) are observable, we
have to characterize them further to make them identifiable. First we
assume that no linear combinations of ${\mathbf x}_t$ are white noise,
as any such components can be absorbed into~$\bve_t$
[see condition (C1) below]. We also assume
that the rank of $\bA$ is $r$. Otherwise~(\ref{c1}) may be expressed
equivalently in terms of a~lower-dimensional factor. Furthermore, since
(\ref{c1}) is unchanged if we replace $(\bA , {\mathbf x}_t)$ by
$(\bA\bH, \bH^{-1}{\mathbf x}_t)$ for any invertible $r\times r$
matrix~$\bH$, we may assume that the columns of $\bA= ( {\mathbf a}_1,
\ldots, {\mathbf a}_r)$ are orthonormal, that is, $\bA' \bA= \bI_r$,
where $\bI_r$ denotes the $r\times r$ identity matrix. Note that even
with this constraint, $\bA$ and ${\mathbf x}_t$ are not uniquely
determined in~(\ref{c1}), as the aforementioned replacement is still
applicable for any orthogonal $\bH$. However, the factor loading space,
that is, the $r$-dimensional linear space spanned by the columns of
$\bA$, denoted by $\calM(\bA)$, is uniquely defined.

We summarize into condition (C1) all the assumptions introduced so far:
\begin{longlist}[(C1)]
\item[(C1)]
In model~(\ref{c1}), $\bve_t \sim \WN(\bmu_{\ve}, \bSigma_{\ve})$. If
${\mathbf c}'{\mathbf X}_t$ is white noise for a constant
\mbox{${\mathbf c} \in \RR^p$}, then ${\mathbf c}'\cov({\mathbf
X}_{t+k}, \bve_t) = \mathbf{0}$ for any nonzero integers $k$.
Furthermore \mbox{$ \bA' \bA = \bI_r$}.
\end{longlist}

The key for the inference for model~(\ref{c1}) is to determine the number
of factors $r$ and to estimate the $p\times r$ factor loading matrix
$\bA$, or more precisely the factor loading space $\calM(\bA)$. Once
we have obtained an estimator,
say,~$\wh\bA$, a natural estimator
for the factor process is
%
\begin{equation} \label{c2}
\wh{\mathbf x}_t = \wh\bA' \by_t,
\end{equation}
and the resulting residuals are
%
\begin{equation} \label{c3}
\wh\bve_t = (\bI_d - \wh\bA\wh\bA') \by_t.
\end{equation}
The dynamic modeling for $\by_t$ is achieved via such a modeling for
$\wh{\mathbf x}_t$
and the relationship $\wh\by_t = \wh
\bA\wh{\mathbf x}_t$. A parsimonious
fitting for $\wh{\mathbf x}_t$
may be obtained by rotating $\wh{\mathbf x}_t$ appropriately~\cite{TT89}.
Such a\vspace*{2pt} rotation is equivalent to replacing $\wh\bA$ by $\wh\bA\bH$
for an appropriate
$r\times r$ orthogonal matrix
$\bH$. Note that $\calM(\wh\bA) = \calM(\wh\bA\bH)$, and the
residuals~(\ref{c3}) are unchanged with
such a replacement.\looseness=1

\subsection{Estimation for $\bA$ and $r$}
\label{sec4}
An innovation expansion algorithm is proposed in~\cite{PY08}
for estimating $\bA$ based on
solving a sequence of nonlinear optimization problems with at most $p$
variables. Although the algorithm
is feasible for small or moderate $p$ only, it can handle the
situations when the
factor process ${\mathbf x}_t $ is nonstationary.
We outline the key idea below, as our computationally more efficient
estimation method for
stationary cases is based on the same principle.

Our goal is to estimate $\calM(\bA)$, or, equivalently, its orthogonal
complement $\calM(\bB)$, where
$\bB= ({\mathbf b}_1, \ldots, {\mathbf b}_{p-r})$ is a $p \times(p-r)$
matrix for which $(\bA, \bB)$ forms a $p\times p$
orthogonal matrix, that is,
$\bB' \bA= 0$ and $\bB'\bB= \bI_{p-r}$ [see also~(C1)].
It follows from~(\ref{c1}) that
%
\begin{equation} \label{d1}
\bB' \by_t = \bB' \bve_t,
\end{equation}
implying that for any $1\le j \le p-r$, $\{ {\mathbf b}_j' \by_t,   t=0,
\pm1, \ldots\}$ is a white-noise process. Hence, we may search for
mutually orthogonal directions ${\mathbf b}_1, {\mathbf b}_2, \ldots$
one by one
such that the projection of $\by_t$ on each of those directions is
a~white noise. We stop the search when such a direction is no longer
available, and take $p-k$ as the estimated value of $r$, where $k$
is the number of directions obtained in the search.
This is essentially how~\cite{PY08}
accomplish the estimation.
It is irrelevant in the above derivation if ${\mathbf x}_t$ is
stationary or not.

However, a much simpler method is available when ${\mathbf x}_t $,
therefore also $\by_t$, is
stationary:
%
\begin{longlist}[(C2)]
\item[(C2)]
${\mathbf x}_t$ is weakly stationary, and Cov$({\mathbf x}_t, \bve
_{t+k})=0$ for any $k \ge0$.
\end{longlist}
In most factor modeling literature, ${\mathbf x}_t$ and $\bve_s$ are
assumed to be uncorrelated for any $t$ and $s$. Condition (C2)
requires only that the future white-noise components are
uncorrelated with the factors up to the present. This enlarges the
model capacity substantially. Put
\begin{eqnarray*}
\bSigma_y(k) &=& \cov( \by_{t+k}, \by_t),\qquad
\bSigma_x(k) = \cov( {\mathbf x}_{t+k}, {\mathbf x}_t),\\
\bSigma_{x\ve}(k) &=& \cov( {\mathbf x}_{t+k}, \bve_t).
\end{eqnarray*}
It follows from~(\ref{c1}) and (C2) that
%
\begin{equation} \label{d2}
\bSigma_y(k) = \bA\bSigma_x(k) \bA' + \bA\bSigma_{x\ve}(k),
\qquad  k \ge1.
\end{equation}
For a prescribed integer $k_0 \ge1$, define
%
\begin{equation} \label{d3}
\bM= \sum_{k=1}^{k_0} \bSigma_y(k) \bSigma_y(k)'.
\end{equation}
Then $\bM$ is a $p\times p$ nonnegative matrix. It follows from~(\ref{d2}) that
$\bM\bB=0$, that is, the columns of $\bB$ are the eigenvectors of
$\bM$ corresponding
to zero-eigenvalues. Hence conditions (C1) and (C2) imply:

\begin{quote}
\normalsize{\textit{The factor loading space $\calM(\bA)$ is spanned by
the eigenvectors of $\bM$ corresponding to
its nonzero eigenvalues, and the number of the nonzero eigenvalues is
$r$.}}\vadjust{\goodbreak}
\end{quote}

\noindent We take the sum in the
definition of $\bM$ to accumulate the information
from different time lags. This is useful especially when the sample
size $n$ is small.
We use the nonnegative
definite matrix $\bSigma_y(k) \bSigma_y(k)'$ [instead of
$\bSigma_y(k)$]
to avoid the cancellation of the information from different lags. This
is guaranteed by the fact that
for any matrix $\bC$, $\bM\bC=0$ if and only if $\bSigma_y(k)'\bC
=0$ for all $1\le k \le k_0$.
We tend to use small $k_0$, as the autocorrelation is often at its
strongest at
the small time lags. On the other hand, adding more terms will not
alter the value of $r$,
although the estimation for $\bSigma_y(k)$ with large~$k$
is less accurate. The simulation results reported in~\cite{LYB11}
also confirm that the estimation for $\bA$ and $r$, defined below, is
not sensitive to the choice of~$k_0$.

To estimate $\calM(\bA)$, we only need to perform an eigenanalysis on
%
\begin{equation} \label{d4}
\wh\bM= \sum_{k=1}^{k_0} \wh\bSigma_y(k) \wh\bSigma_y(k)',
\end{equation}
where $\wh\bSigma_y(k)$ denotes
the sample covariance matrix of $\by_t$ at lag $k$. Then the estimator
$\wh r$ for the number of factors is defined in~(\ref{d5}) below.
The columns of the estimated factor loading matrix
$\wh\bA$ are the $\wh r$ orthonormal eigenvectors of $\wh\bM$
corresponding to its
$\wh r$ largest eigenvalues.
Note that the estimator~$\wh\bA$ is essentially the same as that
defined in Section 2.4 of
\cite{LYB11},
although a canonical form of the model is used there in order
to define the factor loading matrix uniquely.

Due to the random fluctuation in a finite sample, the estimates for the
zero-eigenvalues of $\bM$ are unlikely to be 0 exactly.
A common practice is to plot all the estimated eigenvalues in a
descending order, and look for
a cut-off value $\wh r$ such that the $(\wh r+1)$th largest eigenvalue
is substantially
smaller than the $\wh r$ largest eigenvalues. This is effectively an
eyeball-test.
The ratio-based estimator defined below may be viewed as an enhanced
eyeball-test,
based on the same idea as~\cite{W10}.
In fact this ratio-based estimator benefits from the faster convergence
rates of the
estimators for the zero-eigenvalues; see Proposition~\ref{prop1} in
Section~\ref{sec51} below,
and also Theorems~\ref{thm1} and~\ref{thm2}
in Section~\ref{sec52} below.
The other available methods for determining $r$ include the information
criteria approaches of~\cite{BN02,BN07} and~\cite{HL07}, and
the bootstrap approach
of~\cite{BYZ10},
though the settings considered in
those papers are different.

\textit{A ratio-based estimator for $r$}.
We define an estimator for the number of factors $r$ as follows:
%
\begin{equation}\label{d5}
\wh r = \mathop{\arg\min}_{1\le i \le R} \wh\la_{i+1}/\wh\la_{i},
\end{equation}
where $\wh\la_1 \ge\cdots\ge\wh\la_p$ are the eigenvalues of
$\wh\bM$,
and $ r< R <p$ is a constant.

In practice we may use, for example, $R=p/2$. We cannot extend the
search up to $p$, as the minimum eigenvalue of $\wh\bM$ is likely to
be practically~0, especially when~$n$ is small and~$p$ is large. It
is worthy noting that when~$p$ and~$n$ are on the same order, the
estimators for eigenvalues are no longer consistent. However, the
ratio-based estimator~(\ref{d5}) still works well. See
Theorem~\ref{thm2}(iii) below.

The above estimation methods for $\bA$ and $r$
can be extended to
those nonstationary time series for which a generalized lag-$k$
autocovariance matrix 
is well defined (see, e.g.,~\cite{PP06}).
In fact, the methods are still applicable
when the weak limit of
the generalized lag-$k$ autocovariance matrix
%
\[
\hat\bS_y(k) = n^{-\alpha} \sum_{t=1}^{n-1}(\by_{t+k} - \bar{\by
})(\by_t - \bar{\by})'
\]
exists for $1\le k \le k_0$,
where $\alpha> 1$ is a constant.
Further developments on those lines will be reported elsewhere.
For the factor modeling for high-dimensional volatility processes based
on a similar idea,
we refer to~\cite{PPPY11,TWYZ11}.

\section{Estimation properties}
\label{sec5}

Conventional asymptotic properties are established
under the setting that the sample size $n$ tends to $\infty$ and
everything else remains fixed. Modern time series analysis
encounters the situation when the number of time series $p$ is as
large as, or even larger than, the sample size $n$.
Then the asymptotic properties established under the setting when both
$n$ and $p$ tend to $\infty$ are more relevant. We deal with these two
settings in Section~\ref{sec51} and Sections~\ref{sec52}--\ref{sec525}
separately.

\subsection{\texorpdfstring{Asymptotics when $n \to\infty$ and $p$ fixed}{Asymptotics when n -> infinity and $p$ fixed}}
\label{sec51}

We first consider the asymptotic properties under the assumption
that $n \to\infty$ and $p$ is fixed. These properties reflect the
behavior of our
estimation method in the cases when~$n$ is large and $p$ is small.
We introduce some regularity conditions first. Let $\la_1, \ldots,
\la_p$ be the eigenvalues
of the matrix $\bM$:
\begin{longlist}[(C3)]
\item[(C3)] $\by_t$ is strictly stationary and $\psi$-mixing with the
mixing coefficients~$\psi(\cdot)$ satisfying the condition that
$\sum_{t\ge1} t \psi(t)^{1/2} < \infty$. Furthermore, $E\{ |\by_t|^4
\} < \infty$ element-wisely.
\item[(C4)] $\la_1 > \cdots> \la_r > 0 = \la_{r+1} = \cdots= \la_p$.
\end{longlist}

Section 2.6 of~\cite{FY03}
gives a compact survey on the mixing properties of time series. The use
of the $\psi$-mixing condition in (C3) is for technical convenience.
Note that $\bM$ is a nonnegative
definite matrix. All its eigenvalues are nonnegative. Condition (C4) assumes
that its $r$ nonzero eigenvalues are distinct from each other. While this
condition is not essential, it substantially simplifies the
presentation of the convergence properties
in Proposition~\ref{prop1} below.
Let $\bgamma_j$ be a unit eigenvector of $\bM$ corresponding to the
eigenvalue $\la_j$. We denote by
$(\wh\la_1, \wh\bgamma_1), \ldots, (\wh\la_p, \wh\bgamma_p)$ the
$p$ pairs of eigenvalue and eigenvector of matrix $\wh\bM$: the
eigenvalues $\wh\la_j$ are arranged in descending order, and
the eigenvectors $\wh\bgamma_j$ are orthonormal. Furthermore, it may
go without explicit statement that $\wh\bga_j$ may be replaced by
$-\wh\bga_j$ in order to match the direction of $\bga_j$ for $1\le j
\le r$.
%
\begin{proposition}
\label{prop1}
Let conditions \textup{(C1)--(C4)} hold. Then as $n\to\infty$ (but $p$ fixed), it
holds that:
\begin{longlist}
\item $|\wh\la_j - \la_j | = O_P(n^{-1/2})$ and
$\|\wh\bgamma_j - \bgamma_j \| = O_P(n^{-1/2})$ for $j=1, \ldots,
r$, and

\item $\wh\la_j = O_P(n^{-1})$ for $j= r+1, \ldots, p$.
\end{longlist}
\end{proposition}

The proof of the above proposition is in principle the same as that
of Theorem~1 in~\cite{BYZ10}, and
is therefore
omitted.
\subsection{\texorpdfstring{Asymptotics when $n \to\infty,   p \to\infty$ and $r$ fixed}{Asymptotics when n -> infinity, p -> infinity and r fixed}}
\label{sec52}

\begin{figure}

\includegraphics{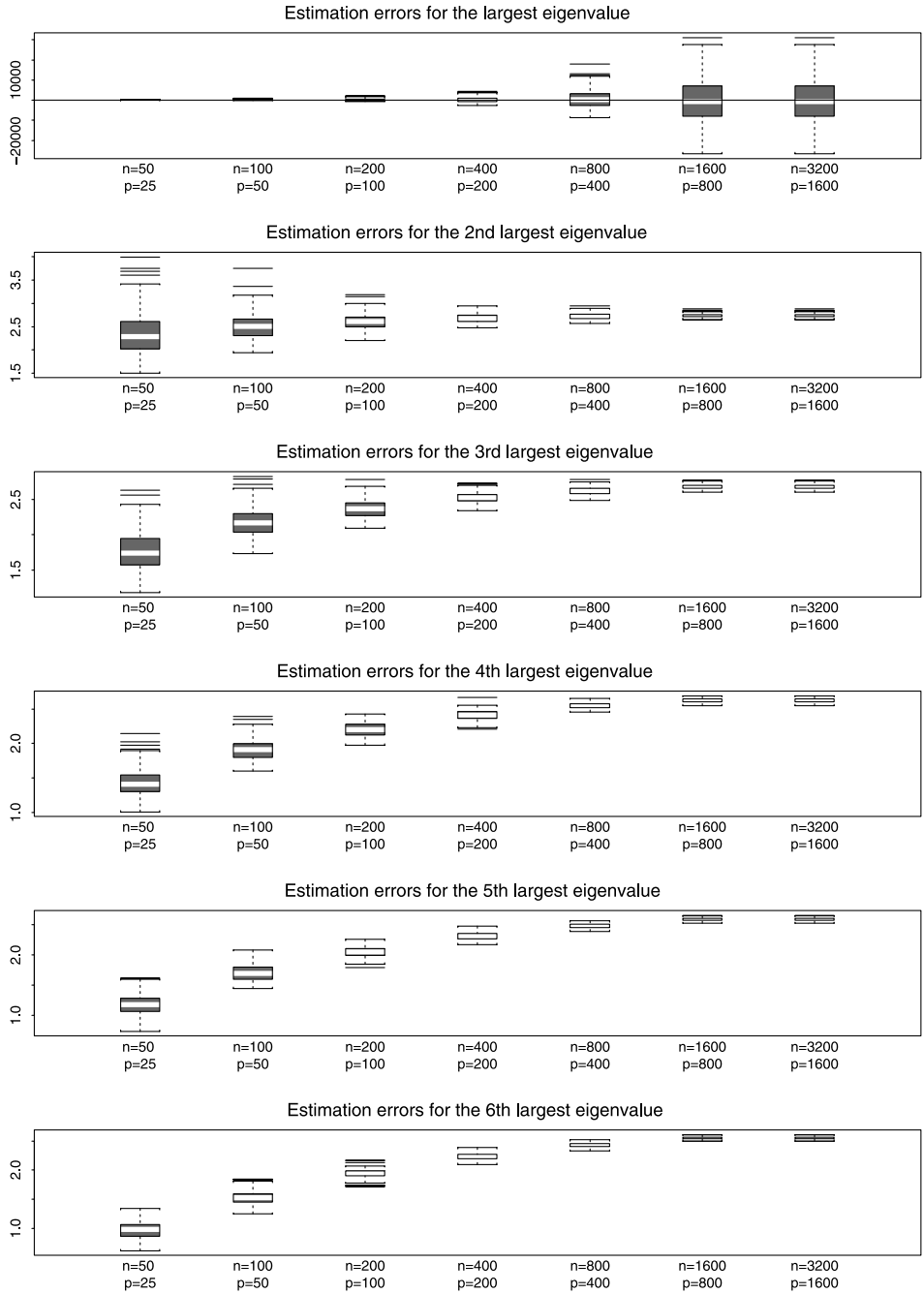}

\caption{Boxplots for the errors in
estimating the first six eigenvalues of $\bM$ with $r=1$ and all the
factor loading
coefficients being 1.}
\label{largeP}
\end{figure}
\begin{figure}

\includegraphics{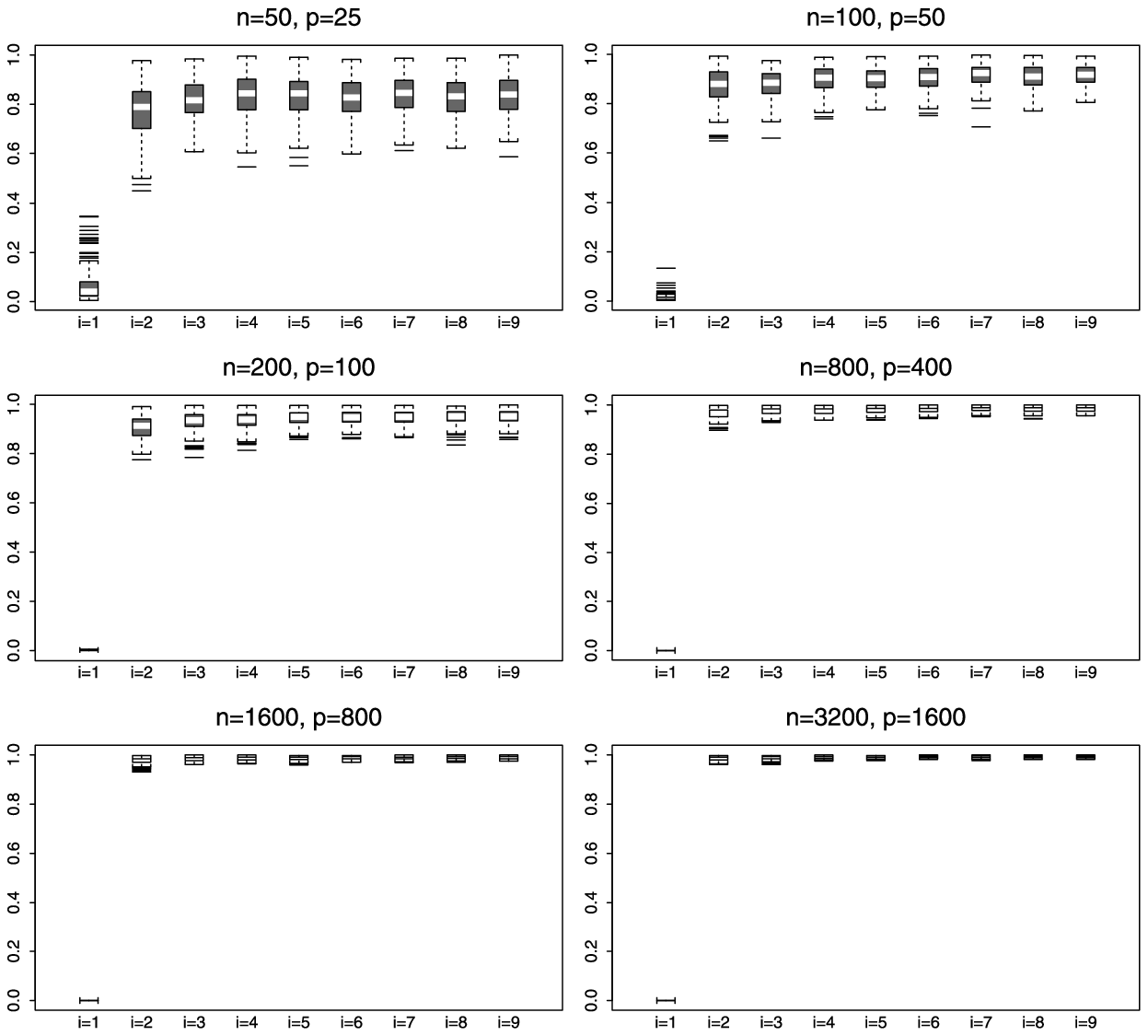}

\caption{Boxplots for the ratios $\wh\la_{i+1}/\wh\la_i$, with
$r=1$ and
all the factor loading coefficients being 1.}
\label{fig5}
\end{figure}

To highlight the radically different behavior when $p$ diverges
together with $n$,
we first conduct some simulations: we set in model~(\ref{c1}) $r=1$,
$\bA'= (1, \ldots, 1)$, $\bve_t$ are independent $N(0, \bI_p)$, and
${\mathbf x}_t = x_t$ is an AR(1) process defined by $ x_{t+1} = 0.7
x_t + e_t. $
We set the sample size $n=50, 100, 200, 400$,
$800, 1600$ and 3200, and the dimension fixed at half the sample size,
that is, $p=n/2$. Let $\bM$ be defined as in~(\ref{d3}) with $k_0 =1$.
For each setting, we draw 200 samples.
The boxplots of the errors $\wh\la_i - \la_i$, $i=1, \ldots, 6$,
are depicted in Figure~\ref{largeP}. Note that $\la_i =0$ for $i\ge
2$, since $r=1$.
%
The figure shows that those estimation errors do not converge to 0. In
fact those errors seem
to increase when $n$ (and also $p=n/2$) increases.
Therefore the classical asymptotic theory (i.e., $n\to\infty$ and $p$
fixed) such as
Proposition~\ref{prop1} above is irrelevant when
$p$ increases together
with $n$. In spite of the lack of consistency in estimating the
eigenvalues, the ratio-based estimator for
the number of factors $r$ ($=$1) defined in~(\ref{d5}) works perfectly
fine for
this example, as shown in Figure~\ref{fig5}.
In fact it is always the case that $\wh r \equiv1$
in all our experiments even when the sample size is as small as
$n=50$; see Figure~\ref{fig5}.

To develop the relevant asymptotic theory,
we introduce some notation first. For any 
matrix $\bG$, let $\|\bG\|$ be the square root of
the maximum eigenvalue of $\bG\bG'$, and $\|\bG\|_{\min}$ be the
square root of the smallest
nonzero
eigenvalue of $\bG\bG'$. We write $a\asymp b$ if $a=O(b)$ and $b=O(a)$.
Recall $\bSigma_x(k)= \cov({\mathbf x}_{t+k}, {\mathbf x}_t)$ and
$\bSigma_{x\epsilon}(k)= \cov({\mathbf x}_{t+k}, \bve_t) $.
Some regularity conditions are now in order:
\begin{longlist}[(C5)]
\item[(C5)] For a constant $\delta\in[0,1]$, it holds that
$\norm{\bSigma_{x}(k)} \asymp p^{1-\delta}\asymp\norm{\bSigma
_{x}(k)}_{\min}$.\vspace*{1pt}
%
\item[(C6)] For $k=0, 1, \ldots, k_0$, $\norm{\bSigma_{x\epsilon
}(k)} = o(p^{1-\delta})$.
\end{longlist}
\begin{Remark}\label{Remark1}
(i) Condition (C5) looks unnatural. It is derived from more natural
conditions~(\ref{e2}) and~(\ref{e3}) below coupled with
the standardization $\bA'\bA= \bI_r$.
Since $\bA= ({\mathbf a}_1, \ldots, {\mathbf a}_r) $ is $p\times r$
and $p \to\infty$ now,
it is natural to let the norm of each column of $\bA$, before
standardizing to $\bA'\bA= \bI_r$, tend to $\infty$ as well.
To this end, we assume that
%
\begin{equation} \label{e2}
\|{\mathbf a}_j\|^2 \asymp p^{1-\delta_j},  \qquad j =1, \ldots, r,
\end{equation}
where $\delta_j \in[0, 1]$ are constants. We take $\delta_j$ as a
measure of the strength of the factor $x_{tj}$.
We call $x_{tj}$ a strong factor
when $\delta_j=0$, and a weak factor when $\delta_j>0$.
Since $r$ is fixed, it is also reasonable to assume that
for $k=0, 1, \ldots, k_0$,
%
\begin{equation} \label{e3}
|\bSigma_x(k)| \ne0.
\end{equation}
Then condition (C5) is entailed by the standardization $\bA'\bA=\bI_r$
under conditions~(\ref{e3})
and~(\ref{e2}) with $\delta_j= \delta$ for all $j$.

(ii) The condition assumed on $\bSigma_{x \epsilon}(k)$ in (C6)
requires that the correlation between ${\mathbf x}_{t+k}$ ($k\ge0$) and
$\bve_t$ is not too strong. In fact under a natural condition that
$\bSigma_{x\epsilon}(k) = O(1)$ element-wisely, it is implied by
(\ref{e2}) and the standardization $\bA'\bA=\bI_r$ [hence now
$x_{j,t} = O_P(p^{(1-\delta)/2})$ as a result of such standardization]
that $
\norm{\bSigma_{x\epsilon}(k)} = O(p^{1-\delta/2})$.

Now we deal with the convergence rates of the estimated eigenvalues, and
establish the results in the same spirit as Proposition~\ref{prop1}.
Of course the convergence (or divergence) rate for each estimator $\wh
\la_i$ is
slower, as the number of estimated parameters goes to infinity now.
\end{Remark}
%
\begin{theorem}
\label{thm1}
Let conditions \textup{(C1)--(C6)} hold and $h_n=
p^\delta n^{-1/2} \to0$. Then as $n\to\infty$ and $p \to\infty$,
it holds that:
\begin{longlist}
\item $|\hat\lambda_i - \lambda_i| = O_P(p^{2-\delta}n^{-1/2})$
for $i=1,\ldots,r$, and\vspace*{2pt}
\item $\hat\lambda_j = O_P(p^2n^{-1})$ for $j = r+1,\ldots, p$.
\end{longlist}
\end{theorem}
%
\begin{corollary}
\label{cor0}
Under the condition of Theorem~\ref{thm1}, it holds that
\[
\hat\lambda_{j+1}/\hat\lambda_j \asymp1    \qquad\mbox{for }
j=1,\ldots,r-1  \quad \mbox{and}\quad   \hat\lambda_{r+1}/\hat\lambda
_r =
O_P(p^{2\de}/n) \pcon0.
\]
\end{corollary}

The proofs of Theorem~\ref{thm1} and Corollary~\ref{cor0} are presented
in the \hyperref[app]{Appendix}. Obviously when $p$ is fixed, Theorem
\ref{thm1} formally reduces to Proposition~\ref{prop1}. Some remarks
are now in order.
\begin{Remark}\label{Remark2}
(i) Corollary~\ref{cor0} implies that the plot of
ratios $\wh\la_{i+1} /\wh\la_i$, $i=1,2, \ldots,$ will drop
sharply at $i=r$.
This provides a partial theoretical underpinning for the estimator $\wh r$
defined in~(\ref{d5}). Especially when all factors are strong (i.e.,
$\de=0$),
$\wh\la_{r+1}/\wh\la_r = O_P(n^{-1})$. This convergence rate is
independent of
$p$, suggesting that the
estimation for $r$ may not suffer as $p$ increases.
In fact when all the factors are strong, the estimation for $r$ may
improve as $p$ increases.
See Remark~\ref{Remark3}(iv) in Section~\ref{sec525} below.

(ii) Unfortunately, we are unable to derive an explicit
asymptotic expression for the ratios $\wh\la_{i+1} /\wh\la_i$ with $i>r$,
although we make the following conjecture:
%
\begin{equation} \label{e4}
\wh\la_{j+1} / \wh\la_j \pcon1,\qquad
j=(k_0+1) r+1, \ldots, (k_0+1) r+K,
\end{equation}
where $k_0$ is the number of lags used in defining matrix $\bM$ in
(\ref{d3}), and $K\ge1$ is any fixed integer. See also
Figure~\ref{fig5}. Further simulation results, not reported
explicitly, also conform with~(\ref{e4}). This conjecture arises
from the following observation: for $j> (k_0+1)r$, the $j$th
largest eigenvalue of~$\wh\bM$ is predominately contributed by the
term $\sum_{k=1}^{k_0} \wh\bSigma_\ve(k) \wh\bSigma_\ve(k)'$ which
has a~cluster of largest eigenvalues on the order of $p^2/n^2$,
where $\hat\bSigma_\varepsilon(k)$ is the sample lag-$k$
autocovariance matrix for $\bve_t$. See also Theorem~\ref{thm2}(iii)
in Section~\ref{sec525} below.


(iii) The errors in estimating eigenvalues are on the order of
$p^{2-\delta}n^{-1/2}$ or $p^2n^{-1}$, and both do not necessarily
converge to 0. However, since
\begin{eqnarray}
{ \hat\lambda_j \over|\hat\lambda_i - \lambda_i| }
=O_P( p^\delta n^{-1/2}) = O_P(h_n) = o_P(1)\nonumber\\
&&\eqntext{\mbox{for  any }
1\le i \le r \mbox{ and  } r< j\le p ,}
\end{eqnarray}
the estimation errors for the zero-eigenvalues is asymptotically of
an order of magnitude smaller than those for the nonzero-eigenvalues.
\end{Remark}

\subsection{Simulation} \label{sec524}

To illustrate the asymptotic properties in Section~\ref{sec52} above, we
report some simulation results. We set
in model~(\ref{c1}) $r=3$, $n=50, 100,
200$, $400, 800, 1600$ and 3200, and $p=0.2n,   0.5n$, $0.8n$ and $1.2n$.
All the $p\times r$ elements of $\bA$ are generated
independently from the uniform distribution on the interval $[-1, 1]$ first,
and we then divide each of them
by~$p^{\delta/2}$ to make all three factors
of the strength~$\delta$; see~(\ref{e2}). We generate
factor~${\mathbf x}_t$ from a $3\times1$ vector-AR(1) process with
independent $N(0,1)$ innovations and the diagonal autoregressive
coefficient matrix with 0.6, $-0.5$ and 0.3 as the main diagonal
elements. We let $\bve_t$ in~(\ref{c1}) consist of independent
$N(0,1)$ components and they are also independent across $t$. We set
$k_0 =1 $ in~(\ref{d3}) and~(\ref{d4}). For each setting, we
replicate the simulation 200 times.

\begin{table}
\caption{Relative frequency estimates for $P(\wh r =r)$ in the
simulation with
200 replications}
\label{table1}
\begin{tabular*}{\tablewidth}{@{\extracolsep{\fill}}lccccd{1.3}d{1.3}cc@{}}
\hline
& $\bolds{n}$ & \textbf{50} & \textbf{100} & \textbf{200}
& \multicolumn{1}{c}{\textbf{400}} & \multicolumn{1}{c}{\textbf{800}}
& \multicolumn{1}{c}{\textbf{1600}} & \textbf{3200}\\
\hline
$\delta=0$
& $p =0.2n$ & 0.165 & 0.680 & 0.940 & 0.995 & 1& 1& 1\\
& $p =0.5n$ & 0.410 & 0.800& 0.980& 1 & 1 & 1 & 1\\
& $p =0.8n$ &0.560 & 0.815 & 0.990 & 1 & 1& 1& 1\\
& $p =1.2n$ &0.590 & 0.820 & 0.990 & 1 & 1& 1& 1\\
[4pt]
$\delta=0.5$
& $p =0.2n$ & 0.075 & 0.155 & 0.270 & 0.570 & 0.980 & 1& 1\\
& $p =0.5n$ & 0.090 & 0.285 & 0.285& 0.820 & 0.960 & 1& 1\\
& $p =0.8n$ & 0.060& 0.180 & 0.490 & 0.745 & 0.970 & 1& 1\\
& $p =1.2n$ & 0.090& 0.180 & 0.310 & 0.760 & 0.915 & 1& 1\\
\hline
\end{tabular*}
\end{table}

Table~\ref{table1} reports the relative frequency estimates for the
probability $P(\wh r =r) = P(\wh r =3)$ with $\delta=0 $ and 0.5.
The estimation performs better when the factors are stronger. Even
when the factors are weak (i.e., $\delta=0.5$), the estimation for $r$ is
very accurate for $n\ge800$.
When the factors are strong (i.e., $\delta=0 $), we observe a
phenomenon coined as ``blessing of
dimensionality'' in the sense that the estimation for
$r$ improves as the dimension $p$ increases.
For example, when the sample size $n=100$, the relative frequencies
for $\wh r = r$ are, respectively, 0.68, 0.8, 0.815 and 0.82 for
$p=20, 50, $ 80 and 120. The improvement is due to the increased
information on $r$ from the added components of $\by_t$ when $p$
increases. When $\delta=0.5$, the columns of $\bA$ are $p$-vectors
with the norm $p^{0.25}$ [see~(\ref{e2})]. Hence we may think that
many elements of $\bA$ are now effectively 0. The increase of the
information on the factors is coupled with the increase of ``noise''
when $p$ increases. Indeed, Table~\ref{table1} shows that when
factors are weak as $\delta=0.5$, the estimation for $r$ does not
necessarily improve as $p$ increases.


\begin{figure}

\includegraphics{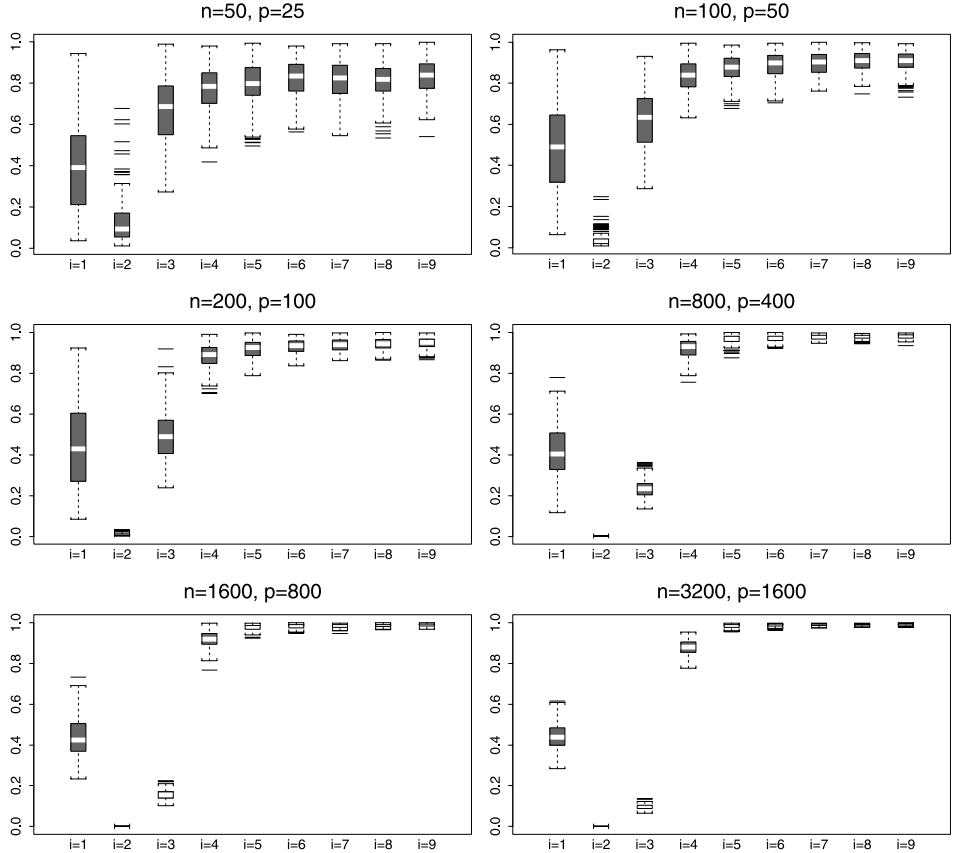}

\caption{Boxplots for the ratios $\wh\la_{i+1}/\wh\la_i$ with
two strong factors ($\delta=0$) and one weak factor ($\delta=0.5$) and
$r=3$, $p=n/2$.}
\label{fig11}
\end{figure}

We also experiment with a setting with two strong factors (with \mbox{$\delta=0$})
and one weak factor (with $\delta=0.5$). Then the ratio-based
estimator $\wh r$
tends to take two values, picking up the two strong factors only.
However Figure~\ref{fig11}
indicates that
the information on the third weak factor is not lost. In fact, $\wh\la
_{i+1}/\wh\la_i$
tends to take the second smallest value at $i=3$. In this case a~two-step
estimation procedure should be employed in order to identify the number
of factors
correctly; see Section~\ref{sec6} below.

\subsection{Improved rates for the estimated eigenvalues}
\label{sec525}

The rates in Theorem~\ref{thm1} can be
further improved,
if we are prepared to entertain some additional conditions on $\bve_t$ in
model~(\ref{c1}). Such an improvement is relevant as the condition
that $h_n = p^\de n^{-1/2} \to0$, required in Theorem~\ref{thm1},
is sometimes
unnecessary. 
For example, in Table~\ref{table1}, the ratio-based estimator $\wh r$
works perfectly well when $\delta=0.5$ and $n$ is sufficiently
large (e.g., $n\ge800$), even though $ h_n = (p/n)^{1/2} \not\to
0$.
Furthermore, in relation to the phenomenon of ``blessing of
dimensionality''
exhibited in Table~\ref{table1}, Theorem~\ref{thm1} fails to reflect
the possible improvement on the estimation for $r$ when $p$ increases;
see also Remark~\ref{Remark2}(i).
We first introduce some additional conditions on $\bve_t$:
\begin{longlist}[(C7)]
\item[(C7)] Let $\varepsilon_{jt}$ denote the $j$th component of $\bve_t$.
Then $\varepsilon_{jt}$ are independent for different $t$ and $j$,
and have mean 0 and common variance $\sigma^2 < \infty$.
\item[(C8)] The distribution of each $\varepsilon_{jt}$ is symmetric.
Furthermore, $E(\varepsilon_{jt}^{2k+1}) = 0$, and $E(\varepsilon_{jt}^{2k})
\leq(\tau k)^k$ for all\vspace*{2pt} $1\le j \le p$ and $t, k \ge1$, where $\tau
>0$ is
a constant independent of $j, t, k $.
\item[(C9)] All the eigenvalues of $\bSigma_{\bve}$
are uniformly bounded as $p \rightarrow\infty$.
\end{longlist}
The moment condition $E(\varepsilon_{jt}^{2k}) \leq(\tau k)^k$ in (C8)
implies that $\varepsilon_{jt}$ are sub-Gaussian.
Condition (C9) imposes some constraint on the correlations among the
components of $\bve_t$.
When all components of $\{\bve_t\}$ are independent $N(0, \sigma^2)$,
(C7)--(C9) hold. See also conditions (i$'$)--(iv$'$) of~\cite{P09}.
%
\begin{theorem}
\label{thm2}
Let conditions \textup{(C1)--(C8)} hold, $\ell_n \equiv p^{\delta/2}n^{-1/2} \to
0$ and
$n = O(p)$.
Then as $p, n \to\infty$,
the following assertions hold:

\begin{longlist}
\item
$|\hat\lambda_j - \lambda_j| = O_P(p^{2-3\delta/2}n^{-1/2}) $
for $j=1,\ldots,r$,\vspace*{1pt}

\item $\hat\lambda_j = O_P(p^{2-\delta}n^{-1}) $
for $j = r+1,\ldots,(k_0+1)r$,

\item $\hat\lambda_j = O_P(p^{2}n^{-2}) $
for $j=(k_0+1)r + 1,\ldots, p$.

If in addition \textup{(C9)} holds, the rate in \textup{(ii)} above can be
further improved to
%
\begin{equation} \label{new}
\hat\lambda_j = O_P(p^{3/2-\delta}n^{-1/2}),\qquad
j = r+1,\ldots,(k_0+1)r.
\end{equation}
\end{longlist}
\end{theorem}
%
\begin{corollary}
\label{cor1}
Under the conditions of Theorem~\ref{thm2}, it holds that
\[
\hat\lambda_{j+1}/\hat\lambda_j \asymp1,\qquad j=1,\ldots,r-1
\mbox{,\quad and }\quad \hat\lambda_{r+1}/\hat\lambda_r = O_P(p^{\delta}n^{-1}).
\]
If in addition \textup{(C9)} also holds, $\hat\lambda_{r+1}/\hat\lambda_r =
O_P(p^{\delta-1/2}n^{-1/2})$.
\end{corollary}

The proofs of Theorem~\ref{thm2} and Corollary~\ref{cor0} are given
in the \hyperref[app]{Appendix}.
\begin{Remark}\label{Remark3}
%
(i) By comparing with Theorem~\ref{thm1}, the error rate for
nonze\-ro~$\la_j$ in Theorem~\ref{thm2} is improved by a factor
$p^{-\de/2}$, the error rate for zero-eigenvalues is by a factor
$p^{-\de}$ at
least. However, those estimators themselves may still
diverge, as illustrated in Figure~\ref{largeP}.

(ii) Theorem~\ref{thm2}(iii) is an interesting consequence of the
random matrix theory. The key message here is as follows: for the
eigenvalues corresponding purely to the matrix $ \sum_{k=1}^{k_0}
\hat\bSigma_\varepsilon(k)\hat\bSigma_\varepsilon(k)' $, their
magnitudes adjusted for $p^{2-2\delta}$ converge at a super-fast rate.
The matrix $ \sum_{k=1}^{k_0}
\hat\bSigma_\varepsilon(k)\hat\bSigma_\varepsilon(k)' $ is a part of
$\hat\bM$ in~(\ref{d4}), where $\hat\bSigma _\varepsilon(k)$ is the
sample lag-$k$ autocovariance matrix for $\{\bve_t\}$. In particular,
when all the factors are strong (i.e., $\de=0$), the convergence rate
is $n^{-2}$. Such a super convergence rate never occurs when $p$ is
fixed.

(iii) Condition $\ell_n \to0$ is mild, and is weaker than condition
$h_n \to0$ required in Theorem~\ref{thm1}.
For example, when $p \asymp n$, this condition is implied by
the condition $\delta\in[0, 1)$.

(iv) With additional condition (C9), $\wh\la_{r+1} / \wh\la_r =
O_P(p^{-1/2} n^{-1/n})$ when
all factors are strong. This shows that the
speed at which $\wh\la_{r+1} / \wh\la_r $ converges to 0 increases when
$p$ increases. This property gives a theoretical explanation why
the identification for $r$ becomes easier for larger $p$ when all
factors are
strong (i.e., $\de=0$). See Table~\ref{table1}.
\end{Remark}

\section{Two-step estimation}
\label{sec6}

In this section, we outline a two-step estimation
procedure. We will show that it is superior than the one-step procedure
presented
in Section~\ref{sec4} for the determination of the number of factors
as well as for the estimation of the factor loading matrices
in the presence of the factors with different degrees of strength.
A similar procedure is described in~\cite{PP06} to improve the
estimation for factor loading matrices in the presence of small eigenvalues,
although they gave no theoretical underpinning on why and when
such a procedure is advantageous.

Consider model~(\ref{c1}) with $r_1$ strong factors with strength
$\delta_1 = 0$ and $r_2$ weak factors with strength $\delta_2 >0$,
where $r_1+r_2=r$. Now~(\ref{c1}) may be written as
%
\begin{equation}\label{e6}
\by_t = \bA{\mathbf x}_t + \bve_t = \bA_1{\mathbf x}_{1t} + \bA
_2{\mathbf x}_{2t} + \bve_t,
\end{equation}
where ${\mathbf x}_t = ({\mathbf x}_{1t}' \enskip  {\mathbf x}_{2t}')'$, $\bA
= (\bA_1 \enskip  \bA_2)$ with
$\bA'\bA= \bI_r$, ${\mathbf x}_{1t}$ consists of $r_1$ strong
factors, and
${\mathbf x}_{2t}$ consists of $r_2$ weak factors. Like model~(\ref{c1}) in Section~\ref{sec31}, $\bA= (\bA_1, \bA_2)$ and ${\mathbf
x}_t = ({\mathbf x}_{1t}, {\mathbf x}_{2t})$ are not uniquely defined,
but only $\calM(\bA)$ is. Hereafter $\bA= (\bA_1, \bA_2)$
corresponds to a suitably rotated version of the original~$\bA$ in
model~(\ref{e6}), where now~$\bA$ contains all the eigenvectors of~$\bM$ corresponding to its nonzero eigenvalues. Refer to~(\ref{d3})
for the definition of~$\bM$.

To present the two-step estimation procedure clearly, let us assume
that we know $r_1$ and
$r_2$ first.
Using the method in Section~\ref{sec4}, we first obtain the estimator
$\wh\bA\equiv(\wh\bA_1, \wh\bA_2)$ for the factor loading matrix
$ \bA= (\bA_1, \bA_2)$,
where the columns of $\wh\bA_1$ are the $r_1$ orthonormal
eigenvectors of
$\wh\bM$ corresponding to its $r_1$ largest eigenvalues.
In practice we may identify $r_1$ using, for example, the
ratio-based estimation method~(\ref{d5}); see Figure~\ref{fig11}. We
carry out the second-step
estimation as follows. Let
%
\begin{equation}\label{e7}
\by_t^* = \by_t - \wh\bA_1\wh\bA_1'\by_t
\end{equation}
for all $t$.
We perform the same estimation for data $\{ \by_t^*\} $ now, and
obtain the $p \times r_2$ estimated factor loading matrix
$\wt\bA_2$ for the $r_2$ weak factors.
Combining the two estimators together, we obtain the final estimator
for $\bA$
as
%
\begin{equation} \label{e8}
\wt\bA= (\wh\bA_1,   \wt\bA_2).
\end{equation}

Theorem~\ref{thm3} below presents the convergence rates for both the
one-step estimator
$\wh\bA= (\wh\bA_1, \wh\bA_2)$
and the two-step estimator $\wt\bA= (\wh\bA_1,   \wt\bA_2)$.
It shows that
$\wt\bA$ converges to $\bA$ at a faster rate than $\wh\bA$.
The results are established with known $r_1$ and $r_2$. In practice
we estimate $r_1$ and $r_2$ using the ratio-based estimators. See also
Theorem~\ref{thm4} below.
We introduce some regularity conditions first.
Let $\bSigma_{12}(k) = \cov({\mathbf x}_{1,t+k},
{\mathbf x}_{2t})$, $\bSigma_{21}(k) =\cov({\mathbf x}_{2,t+k},
{\mathbf x}_{1t})$, $\bSigma_{i}(k) =
\cov({\mathbf x}_{i,t+k}, {\mathbf x}_{it})$ and
$\bSigma_{i\epsilon}(k) = \cov({\mathbf x}_{i,t+k}, \bve_t)$ for $i=1,2$:
\begin{longlist}[(C6)$'$]
\item[(C5)$'$]
For\vspace*{1pt} $i=1,2$, $1 \leq k \leq k_0$,
$\norm{\bSigma_{i}(k)} \asymp p^{1-\delta_i} \asymp\norm{\bSigma
_{i}(k)}_{\min}$,
$\norm{\bSigma_{21}(k)} \asymp\norm{\bSigma_{21}(k)}_{\min}$ and
$\norm{\bSigma_{12}(k)} = O(p^{1-\delta_2/2})$.\vspace*{1pt}


\item[(C6)$'$] $\cov({\mathbf x}_t, \ve_s) =0 $ for any $t,s$.
%
\end{longlist}
The condition on $\bSigma_i(k)$ in (C5)$'$ is an analogue to condition
(C5). See Remark~\ref{Remark1}(i) in Section~\ref{sec52} for the background of those
conditions.
The order of $\norm{\bSigma_{21}(k)}_{\min}$ will
be specified in the theorems below. The order of
$\norm{\bSigma_{12}(k)}$ is not restrictive, since
$p^{1-\delta_2/2}$ is the largest possible order when $\delta_1=0$.
See also the discussion
in Remark~\ref{Remark1}(ii).
Condition (C6)$'$ replaces condition (C6). Here we impose a strong
condition $\bSigma_{i\epsilon}(k) =0$ to highlight the
benefits of the two-step estimation procedure.
See Remark~\ref{Remark4}(iii) below.
Put
\[
\bW_i =
(\bSigma_i(1), \ldots, \bSigma_i(k_0)), \qquad \bW_{21} =
(\bSigma_{21}(1),\ldots,\bSigma_{21}(k_0)).
\]

\begin{theorem}
\label{thm3} Let conditions \textup{(C1)--(C4)},
\textup{(C5)$'$}, \textup{(C6)$'$}, \textup{(C7)} and
\textup{(C8)} hold. Let $n = O(p)$ and
$\kappa_{n} \equiv p^{\delta_2/2}n^{-1/2} \to0$, as $n \to\infty$.
Then it holds that
\[
\norm{\wh\bA_1 - \bA_1} = O_P(n^{-1/2}),\qquad
\norm{\wt\bA_2 - \bA_2} = O_P(\kappa_{n}) = \norm{\wt\bA- \bA}.
\]
Furthermore,
\[
\norm{\wh\bA_2 - \bA_2} = O_P(\nu_{n}) = \norm{\wh\bA- \bA},
\]
if, in addition, $\nu_{n} \to0$ and $p^{c\delta_2}n^{-1/2} \to0$,
where $\nu_n$ and $c$ are defined as follows:
\[
\nu_{n} = \cases{
p^{\delta_2}\kappa_{n}, &\quad
if $\norm{\bW_{21}}_{\min} = o(p^{1-\delta_2})$ $(c=1)$; \vspace*{2pt}\cr
p^{(2c-1)\delta_2}\kappa_{n}, &\quad if $\norm{\bW_{21}}_{\min}
\asymp p^{1-c\delta_2}$ for $1/2 \leq c < 1$,
and\vspace*{2pt}\cr
&\quad $\norm{\bW_1\bW_{21}'} \leq q \norm{\bW_1}_{\min} \norm{\bW_{21}}$
for $0 \leq q <1 $.}
\]
%
\end{theorem}

Note that $\kappa_{n}/\nu_n \to0$.
Theorem~\ref{thm3} indicates that between $\bA_1$ and $\bA_2$, the
latter is more
difficult to estimate, and the convergence rate of an estimator for
$\bA$\vadjust{\goodbreak}
is determined by the rate for $\bA_2$. This is intuitively
understandable as the coefficient vectors for weak factors effectively
contain many zero-components; see~(\ref{e2}). Therefore a
nontrivial proportion of the components of $\by_t$ may
contain little information on weak factors.
When $\norm{\bW_{21}}_{\min} \asymp p^{1-c\delta_2}$, $\norm{\bW
_2}$ is
dominated by $\norm{\bW_{21}}_{\min}$. The condition
$\norm{\bW_1\bW_{21}'} \leq\break q \norm{\bW_1}_{\min} \norm{\bW
_{21}}$ for $0
\leq q < 1$ is imposed to control the behavior of the $(r_1+1)$th to
the $r$th largest eigenvalues of $\bM$ under this situation. If this is
not valid, those eigenvalues can become very small and give a bad
estimator for $\bA_2$, and thus $\bA$. Under this condition, the
structure of the autocovariance for the strong factors, and the structure
of the cross-autocovariance between the strong and weak factors, are
not similar.

Recall that $\la_j$ and $ \wh\la_j$ are the $j$th largest eigenvalue of,
respectively, $\bM$ defined in~(\ref{d3}) and $\wh\bM$ defined in
(\ref{d4}).
We define matrices~$\bM^*$ and $\wh\bM^*$ in the same manner as~$\bM
$ and
$\wh\bM$ but with $\{ \by_t\}$ replaced by
$\{ \by_t^* \}$ defined in~(\ref{e7}), and denote by~$\la_j^*$ and $
\wh
\la_j^*$ the $j$th largest eigenvalue of,
respectively,~$\bM^*$ and $\wh\bM^*$. The following theorem shows
the different behavior of the ratio of eigenvalues under the one-step
and two-step estimation. Readers who are interested in the explicit
rates for the eigenvalues are referred to Lemma~\ref{lemma} in the \hyperref[app]{Appendix}.
%
\begin{theorem}\label{thm4}
Under the same conditions of Theorem~\ref{thm3}, the
following assertions hold:

\begin{longlist}
\item
For $1 \le i < r_1$ or $r_1 < i < r$, $\wh\la_{i+1} / \wh\la_i
\asymp1$. For $1 \le1 < r_2$, \mbox{$ \wh\la^*_{j+1} / \wh\la^*_j
\asymp1$}.
\item
$\wh\la_{r+1} /\wh\la_r \pcon0$ and $\wh\la_{r_1+1} / \wh
\la_{r_1} =o_p(\wh\la_{r+1} /\wh\la_r )$ provided
\[
\delta_2 > 1/(8c-1), \qquad  p^{(1-\delta_2)/2}n^{-1/2} \to0,\qquad
p^{(6c-1/2)\delta_2 -1/2}n^{-1/2} \to\infty.
\]
\item
$\wh\la_{r+1} /\wh\la_r \pcon0$ and $ \wh\la_{r_2 +1}^* /
\wh\la^*_{r_2} =o_p(\wh\la_{r+1} /\wh\la_r )$ provided
$p^{(4c - 3/2)\delta_2-1/2}\times\allowbreak n^{1/2} \rightarrow\infty$.
\end{longlist}
%
%
%
%
%
%
%
%
\end{theorem}

%
%
%
%
%
%
%
%
\begin{Remark}\label{Remark4}
(i) Theorem~\ref{thm4}(i) and (ii) imply
that the
one-step estimation is likely to lead to $\wh r = r_1$. For instance,
when $p \asymp n$,
then Theorem~\ref{thm4}(ii) says that $\wh\la_{r_1+1} / \wh\la
_{r_1}$ has a faster rate of convergence than $\wh\la_{r+1} /\wh\la
_r$ as long as $\delta_2 > 2/5$. Figure~\ref{fig11} shows exactly
this situation.

(ii) Theorem~\ref{thm4}(iii) implies that the two-step estimation is more
capable to identify the additional $r_2$ factors than the one-step
estimation. In particular, if $p \asymp n$, $\wh\la_{r_2 +1}^* / \wh
\la^*_{r_2}$ always has a faster rate of convergence than $\wh\la_{r+1}
/\wh\la_r$. Unfortunately we are
unable to establish the asymptotic properties for $\wh\la_{i+1}/
\wh\la_i$ for $i>r$, and $\wh\la_{j+1}^*/ \wh\la_j^*$ for
$j>r_2$, though we believe that conjectures similar to~(\ref{e4})
continue to hold.

(iii) When $\delta_1 > 0$ and/or the cross-autocovariances between different
factors and the noise are stronger, the similar and more complex results
can be established via more involved algebra in the proofs.
\end{Remark}



\section{Real data examples} \label{sec7}
We illustrate our method using two real data sets.\vspace*{-2pt}

\begin{Example}\label{Example1}
We first analyze the daily returns of 123 stocks in
the period 2 January 2002--11 July 2008. Those stocks were
selected among those included in the S\&P500 and were traded
every day during the period. The returns were calculated in
percentages based on the daily close prices. We have in total
$n=1642$ observations with $p=123$. We apply the eigenanalysis to
the matrix $\wh\bM$ defined in~(\ref{d4}) with $k_0 =5$. The
obtained eigenvalues (in descending order) and their ratios are
plotted in Figure~\ref{figSPeig}. It is clear that the ratio-based
estimator~(\ref{d5}) leads to $\wh r =2$, indicating two factors.
%
\begin{figure}
\begin{tabular}{@{}c@{}}

\includegraphics{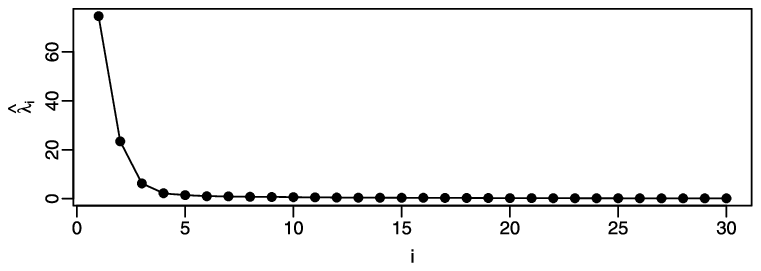}
\\
(a)\\[4pt]

\includegraphics{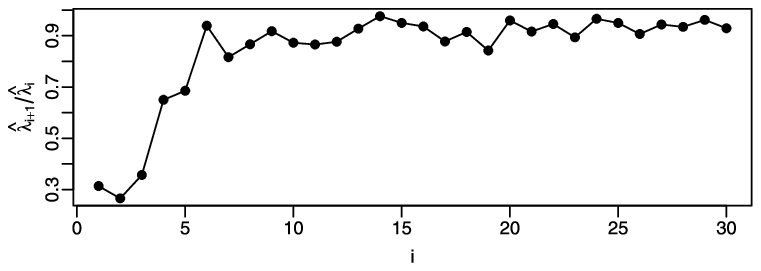}
\\
(b)
\end{tabular}
\caption{Plots of the estimated eigenvalues \textup{(a)} and the
ratios of estimated eigenvalues of~$\wh\bM$ \textup{(b)}~for Example \protect\ref{Example1}.}
\label{figSPeig}
\vspace*{-3pt}
\end{figure}
Varying the value of $k_0$ between 1 and 100 in the definition of
$\wh\bM$ leads to little change in the ratios $\wh\la_{i+1}/\wh
\la_i$, and the estimate $\wh r =2$ remains unchanged.
Figure~\ref{figSPeig}(a) shows that $\wh
\la_{i}$ is close to 0 for all $i\ge5$. Figure~\ref{figSPeig}(b)
indicates that the ratio $\wh\la_{i+1}/\wh\la_i$ is close to 1 for
all large $i$, which is in line with conjecture~(\ref{e4}).

\begin{figure}

\includegraphics{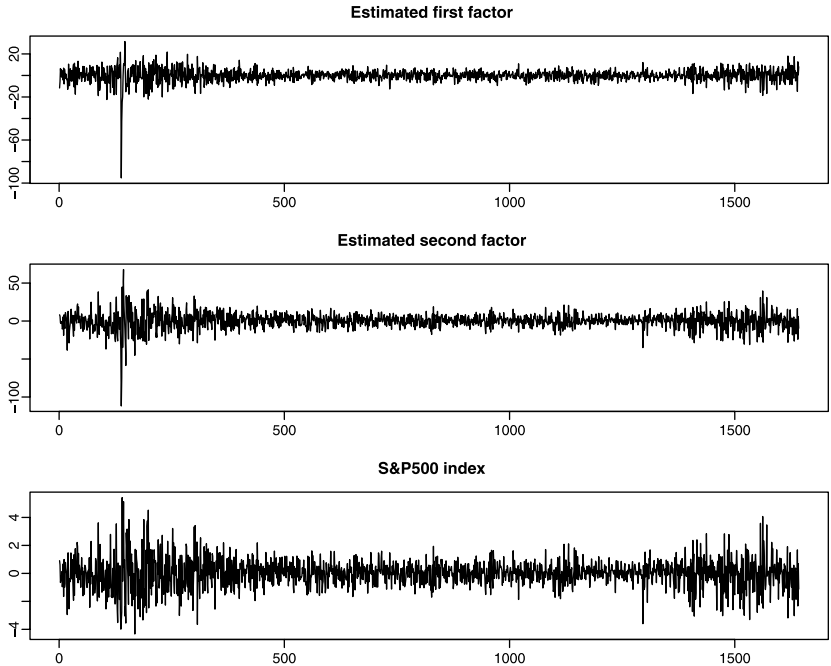}

\caption{The time series plots of the two estimated factors and the
return series
of the S\&P500 index in the same time period.}
\label{figSPTS}
\end{figure}

The first two panels of Figure~\ref{figSPTS} display the time series
plots of the two component series of the estimated factors $\wh
{\mathbf x}_t$ defined as in~(\ref{c2}). Their cross-autocorrelations are
presented in Figure~\ref{figSPfacACF}. Although each of the two
estimated factors shows little significant autocorrelation, there
are some significant cross-correlations between the two series.
The cross-autocorrelations of the three residual series $\wh\bgamma_j'
\by_t$ for $j=3, 4, 5$ are not\vspace*{-1pt} significantly different
from 0, where $\wh\bgamma_j$ is the unit eigenvector of $\wh\bM$
corresponding to its $j$th largest eigenvalue. If there were any serial
correlations left in the data after extracting the two estimated
factors, those correlations are most likely to show up in those three
residual series.\vadjust{\goodbreak}

Figure~\ref{figSPeig} may suggest the existence of a third and weaker
factor, though there are hardly\vspace*{1pt} any significant
autocorrelations in the series $\wh\bgamma_3' \by_t$. In fact $
\wh\la_3 = 6.231$ and $\wh\la_4/\wh\la_3= 0.357$. Note\vspace*{-1pt}
that now $\wh\la_j$ is not necessarily a consistent estimator for
$\la_j$ although $\wh\la_{r+1}/\wh\la_r \pcon0$; see Theorem
\ref{thm1}(ii) and Corollary~\ref{cor0}. To investigate this further,
we apply the two-step estimation procedure presented in
Section~\ref{sec6}. By subtracting the two estimated factors from the
above, we obtain the new data $\by_t^*$ [see~(\ref{e8})]. We
then\vspace*{1pt} calculate the eigenvalues and their ratios of the
matrix $\wh \bM^*$. The minimum value of the ratios is
$\wh\la_2^*/\wh\la_1^* = 0.667$, which is closely followed by
$\wh\la_3^*/\wh\la_2^* = 0.679$ and $\wh\la^*_4/\wh\la_3^* = 0.744$.
There is no evidence to suggest that \mbox{$\wh\la_2^*/\wh\la_1^*
\to0$};
see Theorem~\ref{thm4}. This reinforces our
choice $\wh r =2$.

\begin{figure}

\includegraphics{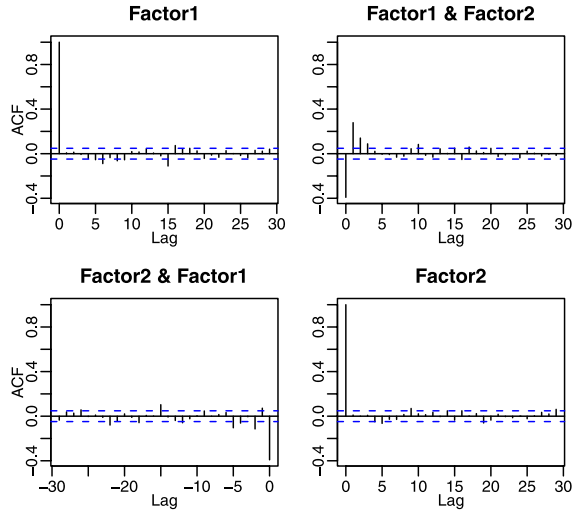}

\caption{The cross-autocorrelations of the two estimated
factors for Example \protect\ref{Example1}.} \label{figSPfacACF}
\end{figure}

With $p$ as large as 123, it is difficult to gain insightful
interpretation on the
estimated factors by looking through the coefficients in $\wh\bA$
[see~(\ref{c2})].
To link our fitted factor model with some classical asset pricing
theory in finance,
we wonder if the market index (i.e., the S\&P500 index) is a factor in
our fitted model,
or more precisely, if it can be written as a linear combination of the
two estimated
factors.
When this is true, $\bP\bu=0$,
where $\bu$ is the $1642\times1$ vector consisting of the returns of
the S\&P500 index over the same
time period, and $\bP$ denotes the projection matrix onto the orthogonal
complement of the linear space spanned by the two component
series $\wh{\mathbf x}_t$, which is a 1640-dimensional subspace in $R^{1642}$.
This S\&P500 return series is plotted together with the two component
series $\wh{\mathbf x}_t$
in Figure~\ref{figSPTS}.
It turns out that $\|\bP\bu\|^2 $ is not exactly 0 but $\|\bP\bu
\|^2/ \|\bu\|^2=0.023$, that is, the 97.7\% of the S\&P500
returns can be expressed as a linear
combination of the two estimated factors.
Thus our analysis suggests the following model for $\by_t$---the daily
returns of the 123 stocks:
\[
\by_t = {\mathbf a}_1 u_t + {\mathbf a}_2 v_t + \bve_t,
\]
where $u_t$ denotes the return of the S\&P500 on the day $t$, $v_t$
is another factor, and $\bve_t$ is a $123\times1$ vector white-noise process.

Figure~\ref{figSPTS} shows that there is an early period with big
sparks in the two estimated
factor processes. Those sparks occurred around 24 September 2002 when
the markets
were highly volatile and the Dow Jones Industrial Average had lost 27\%
of the value
it held on 1 January 2001. However, those sparks are significantly less
extreme in
the returns of the S\&P500 index; see the third panel in Figure~\ref{figSPTS}.
In fact the projected S\&P500 return $\bP\bu$ is the linear
combination of those two
estimated factors 
with the coefficients ($-0.0548, 0.0808$). Two observations may be
drawn from the opposite signs of
those two coefficients: (i) there is an indication that those two factors
draw the energy from the markets with opposite directions, and (ii) the
portfolio S\&P500
index hedges the risks across different markets.
\end{Example}
\begin{Example}\label{Example2}
We analyze a set of monthly average sea surface air pressure records
(in Pascal)
from January 1958 to December 2001
(i.e., 528 months in total)\vadjust{\goodbreak}
over a $10 \times44$ grid in a range of $22.5^\circ$--$110^\circ$ longitude
in the North Atlantic Ocean.
Let $P_t(u,v)$ denote the air pressure in the $t$th month at the location
$(u,v)$, where $ u=1,\ldots,10,  v=1,\ldots,44$ and $t =
1,\ldots,528$.
We first subtract each data point by the
monthly mean over the 44 years at its location:
$\frac{1}{44}\sum_{i=1}^{44}P_{12(i-1) + j}(u,v)$, where $j=1, \ldots,
12$, representing the 12 different months over a year.
We then line up the new data over $10 \times44 =440$ grid points as a
vector $\by_t$, so that
$\by_t$ is a $p$-variate time series with $p=440$. We have $n=528$
observations.

To fit the factor model~(\ref{c1}) to $\by_t$, we calculate the
eigenvalues and
the eigenvectors of the matrix $\wh\bM$ defined in~(\ref{d4}) with
$k_0 =5$. Let $\wh\la_1 > \wh\la_2 > \cdots$ denote the
eigenvalues of $\wh\bM$.
The ratios $\wh\la_{i+1} / \wh\la_i$ are plotted against $i$ in the
top panel of
Figure~\ref{fig13} which indicates the ratio-based estimate for the
number of factor
is $\wh r =1$; see~(\ref{d5}).
However, the second smallest ratio is $\wh\la_4 /\wh\la_3$.
This suggests that there may exist two weaker factors in addition; see
Theorem~\ref{thm4}(ii) and also Figure~\ref{fig11}.
We adopt the two-step estimation procedure presented in Section~\ref{sec6}
to identify the factors of different strength.
By removing the factor corresponding to the largest eigenvalue of $\wh
\bM$,
%
\begin{figure}

\includegraphics{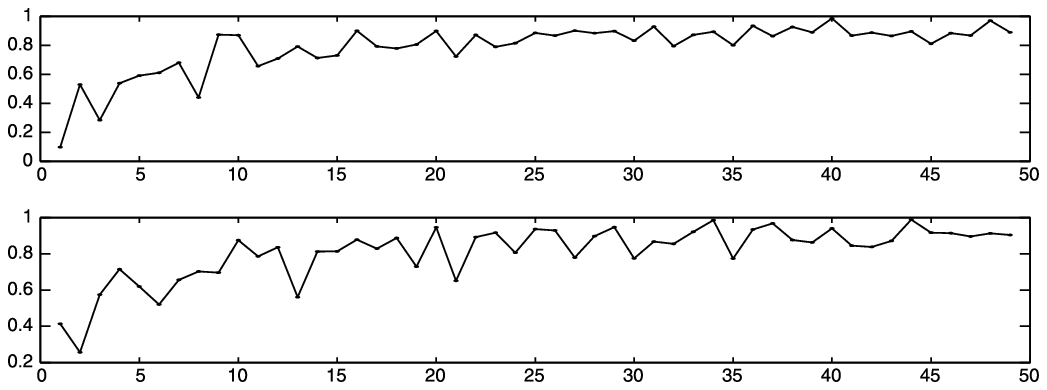}

\caption{Plots of $\wh
\la_{i+1}/\wh\la_i$---the ratio of the eigenvalues of $\wh\bM$
(the top panel) and
$\wh\bM^*$ (the bottom panel), against $i$, for Example \protect\ref{Example2}.}
\label{fig13}
\end{figure}
the resulting ``residuals'' are denoted as $\by_t^*$; see~(\ref{e7}).
Now we repeat the factor modeling for data $\by_t^*$, and plot the
ratios of eigenvalues of matrix $\wh\bM^*$
in the second panel of Figure~\ref{fig13}. It shows clearly the
minimum value at 2,
indicating further two (weaker) factors.
Combining the above two steps together, we set $\hat r = 3$ in the
fitted model.
We repeated the above calculation with $k_0=1$ in~(\ref{d4}). We still
find three factors with the two-step procedure, and the estimated
factors series are very similar to the case when $k_0=5$. This is
consistent with the simulation results in~\cite{LYB11},
where they showed empirically that the estimated factor models are not
sensitive to the choice of $k_0$.

%
\begin{figure}

\includegraphics{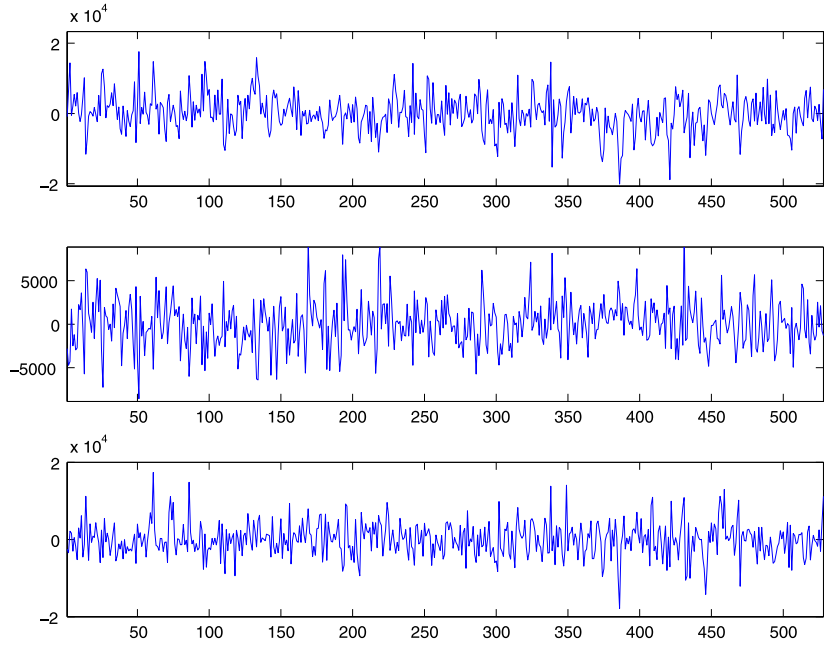}

\caption{Time series plot of the three estimated factors for
Example \protect\ref{Example2}.}\vspace*{-0.5pt}
\label{fig14}
\end{figure}

We present the time series plots for the three estimated factors $\wh
{\mathbf x}_t = \wt\bA' \by_t$ in Figure~\ref{fig14}, where $\wt\bA$ is
a $440\times3$ matrix with the first column being the unit eigenvector
of $\wh\bM$ corresponding to its\vspace*{1pt} largest eigenvalue, and the other two
columns being the orthonormal eigenvectors of $\wh\bM^*$ corresponding
%
\begin{figure}[b]\vspace*{-0.5pt}

\includegraphics{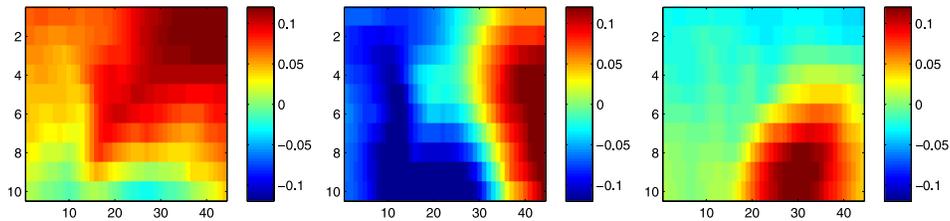}

\caption{Factor loading surface of the first, second and third factors
(from left to right) for Example~\protect\ref{Example2}.}
\label{fig15}
\end{figure}
to its two largest\vadjust{\goodbreak} eigenvalues; see~(\ref{e8}) and also
(\ref{c2}). They collectively account for 85.3\% of the total variation
in $\by_t$ which has $440$ components. In fact each of the three
factors accounts for, respectively, 57.5\%, 18.2\% and 9.7\% of the
total variation of $\by_t$.
Figure~\ref{fig15} depicts the factor loading surfaces of the three
factors. Some interesting regional patterns are observed from those
plots. For example, the first factor is the main driving force for
the dynamics in the north and especially the northeast. The second
factor influences the dynamics in the east and the west in the
opposite directions, and has little impact in the narrow void
between them. The third factor impacts mainly the dynamics of the
southeast region. We also notice that none of those factors can be
seen as idiosyncratic components as each of them affects quite a
large number of locations.\vadjust{\goodbreak}

\begin{figure}

\includegraphics{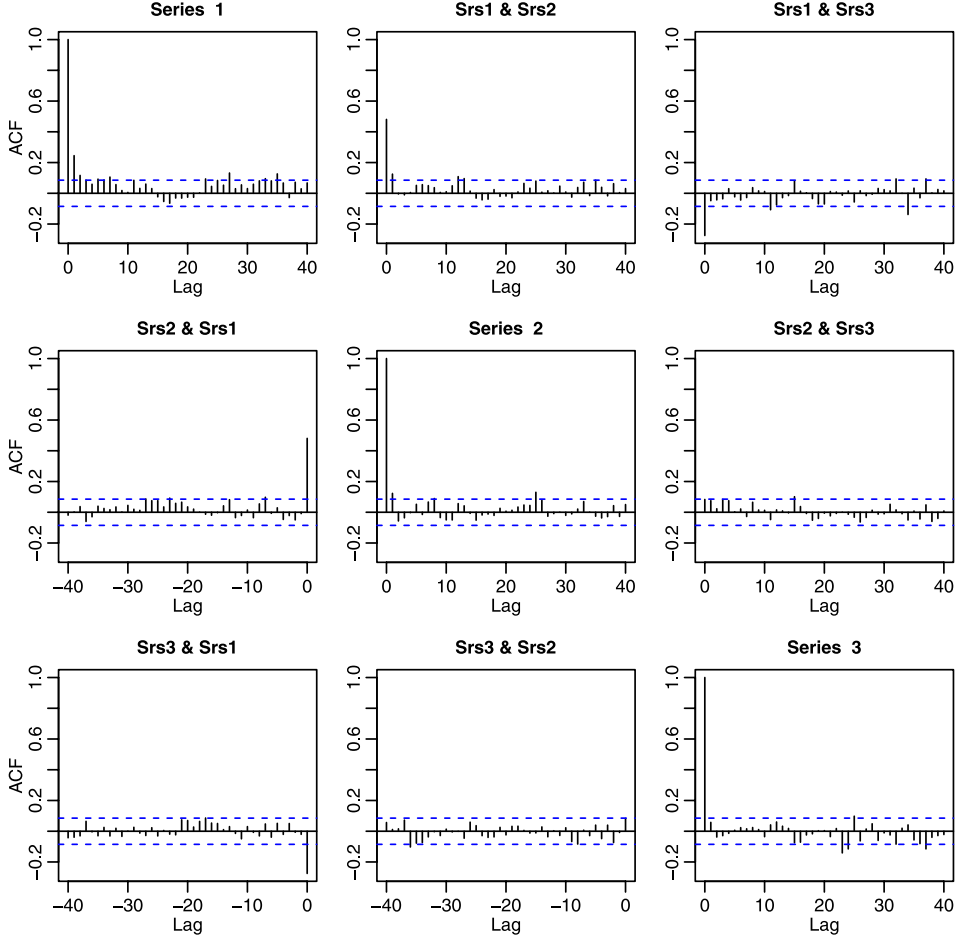}

\caption{Example \protect\ref{Example2}: sample cross-correlation functions for the three
estimated factors.}
\label{fig16}
\end{figure}

\begin{figure}

\includegraphics{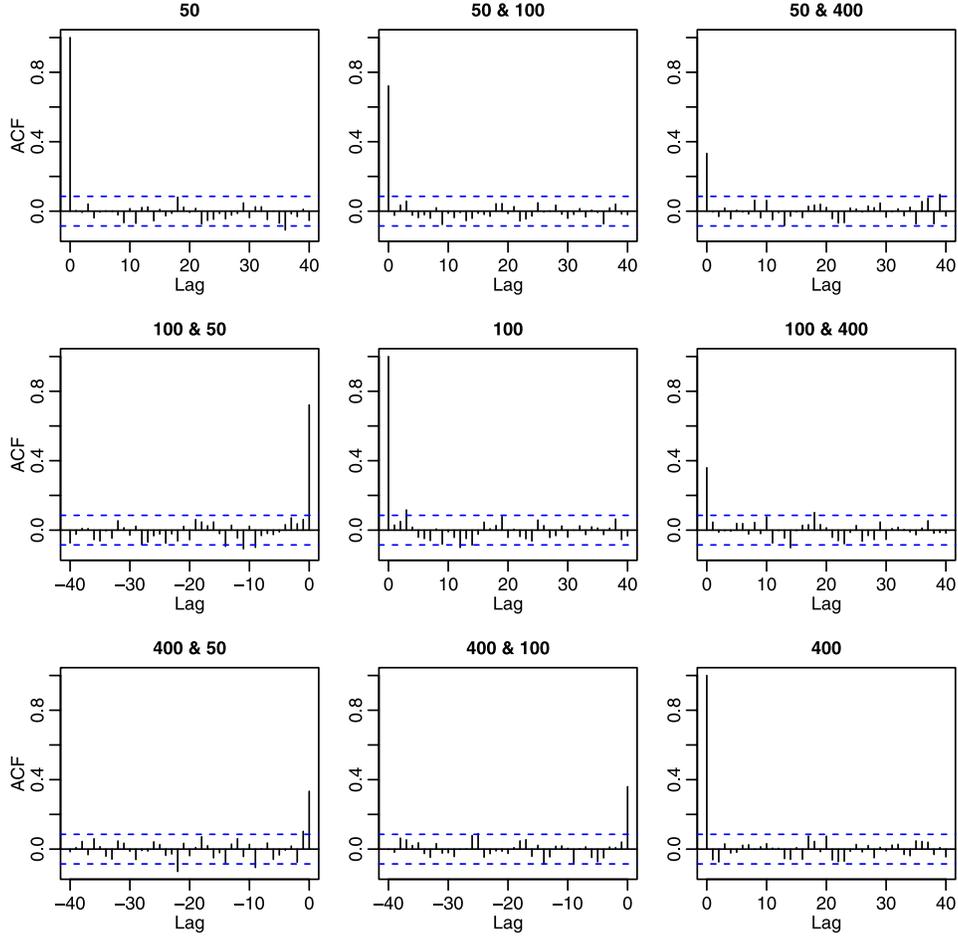}

\caption{Example \protect\ref{Example2}: sample cross-correlation functions for three
residual series. 50 represents grid position (10, 5), 100 for (10, 10)
and 400 for (10, 40).}
\label{fig17}
\end{figure}

Figure~\ref{fig16} presents the
sample cross-correlations for the three estimated factors. It shows
significant, though small, autocorrelations or cross-correla\-tions at
some nonzero lags.
Figure~\ref{fig17} is the sample cross-correlations for three
residuals series
selected from three locations for which one is far apart from the other
two spatially,
showing little autocorrelations at nonzero lags.
This indicates that our approach is capable to identify the factors
based on serial correlations.
\end{Example}

Finally we note that the BIC method of~\cite{BN02}
yields the
estimate $\hat r = n = 528$ for this particular data set. We suspect
that this may be due to the fact that~\cite{BN02}
requires all
the eigenvalues of $\bSigma_\ve$ be uniformly bounded when $p \to
\infty$. This may not be the case for this particular data set, as
the nearby locations are strongly spatially correlated, which may
lead to very large and also very small eigenvalues for
$\bSigma_\ve$. Indeed, for this data set, the three largest
eigenvalues of $\hat\bSigma_\ve$ are on the order of $10^6$, and the
three smallest eigenvalues are practically 0.
Since the typical magnitude of $\hat\bve_t$ is $10^2$ from our
analysis, we have done simulations (not shown here) showing that the
typical largest eigenvalues for $\hat\bSigma_{\ve}$, if $\{\bve_t\}$
is weakly correlated white noise, should be around $10^4$ to $10^5$,
and the smallest around $10^2$ to $10^3$ when $p=440$ and $n=528$.
Such a huge difference in the magnitude of the eigenvalues suggests
strongly that the components of the white-noise vector $\bve_t$ are
strongly correlated. Our method does not require the uniform
boundedness of the eigenvalues of~$\bSigma_{\ve}$.


\begin{appendix}\label{app}
\section*{Appendix}

\begin{pf*}{Proof of Theorem~\ref{thm1}}
We present some notational definitions first. We denote by
$\hat\lambda_j, \hat\bgamma_j$ the $j$th largest eigenvalue of
$\hat\bM$ and the corresponding orthonormal eigenvector, respectively.
The corresponding population values are denoted by $\lambda_j$ and
${\mathbf a}_j$ for the matrix $\bM$. Hence $\hat\bA=
(\hat\bgamma_1,\ldots,\hat\bgamma_r)$ and $\bA= ({\mathbf
a}_1,\ldots,{\mathbf a}_r)$. We also have
\[
\lambda_j = {\mathbf a}_j'\bM{\mathbf a}_j, \qquad\hat\lambda_j =
\hat\bgamma_j'\hat\bM\hat\bgamma_j,\qquad j=1,\ldots,p.
\]

We show some intermediate results now. With conditions (C3) and (C5) and
the fact that $\{\bve_t\}$ is white noise, we have
%
\begin{eqnarray}\label{app16}
\norm{\wh\bSigma_x(k) - \bSigma_x(k)} &=& O_P(p^{1-\delta
}n^{-1/2}), \nonumber\\[-8pt]\\[-8pt]
\norm{\wh\bSigma_{x\epsilon}(k) - \bSigma_{x\epsilon}(k)},\qquad
\norm{\wh\bSigma_{\epsilon x}(k)} &=& O_P(p^{1-\delta/2}n^{-1/2}),
\nonumber
\end{eqnarray}
where $k=0,1,\ldots,k_0$. Then following the proof of Theorem 1 of
\cite{LYB11},
we have the following for $k = 1,\ldots,k_0$:
%
\begin{eqnarray}\label{app1}\quad
&&\norm{\hat\bM- \bM} = O_P\bigl(\norm{\bSigma_y(k)} \cdot\norm{\wh
\bSigma_y(k) - \bSigma_y(k)}\bigr)
\nonumber\\
&&\qquad\mbox{where }\hspace*{40.5pt}
\norm{\bSigma_y(k)} = O(p^{1-\delta}),\nonumber\\[-8pt]\\[-8pt]
&&\qquad\hphantom{\mbox{where }}\norm{\wh\bSigma_y(k) - \bSigma_y(k)} = O_P\bigl(p^{1-\delta}n^{-1/2}
+ p^{1-\delta/2}n^{-1/2} + \norm{\wh\bSigma_\epsilon(k)}\bigr)\nonumber\\
&&\qquad\hphantom{\mbox{where }\norm{\wh\bSigma_y(k) - \bSigma_y(k)}}= O_P\bigl(p^{1-\delta/2}n^{-1/2} +
\norm{\wh\bSigma_\epsilon(k)}\bigr).\nonumber
\end{eqnarray}
Now $\norm{\wh\bSigma_\epsilon(k)} \leq\norm{\wh\bSigma
_\epsilon(k)}_F = O_P(pn^{-1/2})$, where $\norm{\bM}_F =
\operatorname{trace}(\bM\bM')$ denotes the Frobenius norm of $\bM$. Hence from
(\ref{app1}),
%
\begin{eqnarray}\label{app2}
\norm{\wh\bSigma_y(k) - \bSigma_y(k)} &=& O_P(pn^{-1/2})\quad
\mbox{and }\nonumber\\[-8pt]\\[-8pt]
\norm{\hat\bM- \bM} &=& O_P(p^{1-\delta} \cdot pn^{-1/2}) =
O_P(p^{2-\delta}n^{-1/2}).
\nonumber
\end{eqnarray}

For the main proof, consider for $j=1,\ldots,r$, the decomposition
%
\begin{eqnarray}\label{appd1}\quad
&&
\hat\lambda_j - \lambda_j = \hat\bgamma_j'\hat\bM\hat\bgamma
_j -
{\mathbf a}_j'\bM{\mathbf a}_j
= I_1 + I_2 + I_3 + I_4 + I_5\nonumber\\
&&\qquad
\mbox{where }
I_1 = (\hat\bgamma_j - {\mathbf a}_j)'(\hat\bM-
\bM)\hat\bgamma_j,\qquad
I_2 = (\hat\bgamma_j - {\mathbf a}_j)'\bM(\hat\bgamma_j -
{\mathbf a}_j),\nonumber\\[-8pt]\\[-8pt]
&&\qquad\hphantom{\mbox{where }} I_3 = (\hat\bgamma_j - {\mathbf a}_j)'\bM{\mathbf a}_j,\qquad
I_4 = {\mathbf a}_j'(\hat\bM- \bM)\hat\bgamma_j, \nonumber\\
&&\qquad\hphantom{\mbox{where }} I_5 =
{\mathbf a}_j'\bM(\hat\bgamma_j - {\mathbf a}_j).\nonumber
\end{eqnarray}
For $ j=1,\ldots,r$, since $\norm{\hat\bgamma_j - {\mathbf a}_j}
\leq
\norm{\hat\bA- \bA} = O_P({h_n})$ where $h_n = p^{\delta}n^{-1/2}$,
and $\norm{\bM} \leq\sum_{k=1}^{k_0} \norm{\bSigma_y(k)}^2
=O_P(p^{2-2\delta})$ by~(\ref{app1}), together with~(\ref{app2})
we have that
\[
\norm{I_1},   \norm{I_2} = O_P(p^{2-2\delta}h_n^2),  \qquad
\norm{I_3}, \norm{I_4}, \norm{I_5} = O_P(p^{2-2\delta}h_n),
\]
so that
$|\hat\lambda_j - \lambda_j| = O_P(p^{2-2\delta}h_n) =
O_P(p^{2-\delta}n^{-1/2})$, which proves Theorem~\ref{thm1}(i).\vadjust{\goodbreak}

Now consider $j=r+1,\ldots,p$. Define
\[
\wt\bM= \sum_{k=1}^{k_0}\wh\bSigma_y(k)\bSigma_y(k)',\qquad
\hat\bB= (\hat\bgamma_{r+1},\ldots,\hat\bgamma_p),  \qquad
\bB=({\mathbf a}_{r+1},\ldots,{\mathbf a}_p).
\]
Following the same proof of Theorem 1 of~\cite{LYB11},
we can actually show that $\norm{\wh\bB- \bB} = O_P(h_n)$,
so that $\norm{\hat\bgamma_j - {\mathbf a}_j} \leq\norm{\wh\bB-
\bB} =
O_P(h_n)$.

Noting $\lambda_j = 0$ for $j=r+1,\ldots,p$, consider the
decomposition
%
\begin{eqnarray}\label{appd2}
&&\hat\lambda_j = \hat\bgamma_j'\wh\bM\hat\bgamma_j = K_1 + K_2
+ K_3 \nonumber\\
&&\qquad\mbox{where }K_1 = \hat\bgamma_j'(\wh\bM- \wt\bM- \wt\bM' +
\bM)\hat\bgamma_j, \nonumber\\[-8pt]\\[-8pt]
&&\hphantom{\qquad\mbox{where }}K_2 = 2\hat\bgamma_j'(\wt\bM- \bM)(\hat
\bgamma_j -
{\mathbf a}_j),\nonumber\\
&&\hphantom{\qquad\mbox{where }}K_3 =
(\hat\bgamma_j - {\mathbf a}_j)'\bM(\hat\bgamma_j -
{\mathbf a}_j).
\nonumber
\end{eqnarray}
Using~(\ref{app2}),
\[
K_1 = \sum_{k=1}^{k_0} \bigl\|\bigl(\wh\bSigma_y(k) - \bSigma_y(k)\bigr)'\hat
\bgamma_j\bigr\|^2
\leq\sum_{k=1}^{k_0} \norm{\wh\bSigma_y(k) - \bSigma_y(k)}^2 =
O_P(p^2n^{-1}).
\]
Similarly, using~(\ref{app1}) and~(\ref{app2}), and $\norm{\wh\bB-
\bB} = O_P(h_n)$, we can show that
\begin{eqnarray*}
|K_2| &=& O_P(\norm{\wt\bM- \bM} \cdot\norm{\hat\bgamma_j -
{\mathbf a}_j})
= O_P(\norm{\wh\bM- \bM}\cdot\norm{\wh\bB- \bB}) =
O_P(p^2n^{-1}),\\
|K_3| &=& O_P(\norm{\wh\bB- \bB}^2 \cdot\norm{\bM}) =
O_P(p^{2-2\delta}h_n^2) = O_P(p^2n^{-1}).
\end{eqnarray*}
Hence $\hat\lambda_j = O_P(p^2n^{-1})$, and the proof of the theorem is
completed.
\end{pf*}
\begin{pf*}{Proof of Corollary~\ref{cor0}}
The proof of Theorem 1 of~\cite{LYB11}
has shown that (in the notation of this paper)
\[
p^{2-2\delta} = O(\lambda_r).
\]
But we also have
\[
\lambda_r \leq\lambda_1 = \norm{\bM} \leq\sum_{k=1}^{k_0}
\norm{\bSigma_y(k)}^2 = O(p^{2-2\delta}),
\]
where the last equality sign follows from $\norm{\bSigma_y(k)} =
O(p^{1-\delta})$ in~(\ref{app1}). Hence we have $\lambda_i \asymp
p^{2-2\delta}$ for $i=1,\ldots,r$.

Letting $e_i = |\hat\lambda_i - \lambda_i|$ for $i=1,\ldots,r$, we
then have $e_i = O_P(p^{2-\delta}n^{-1/2})$ from Theorem~\ref{thm1}(i). But since $h_n = p^{\delta}n^{-1/2} = o(1)$ implying that
$p^{2-\delta}n^{-1/2} = p^{2-2\delta}h_n = o(p^{2-2\delta})$, we
have $e_i = o_P(\lambda_i)$. Hence we must have $\hat\lambda_i
\asymp\lambda_i \asymp p^{2-2\delta}$ for $i=1,\ldots,r$. This
implies that $\hat\lambda_{j+1}/\hat\lambda_j \asymp1$ for
$j=1,\ldots,r-1$, and together with Theorem~\ref{thm1}(ii),
\[
\hat\lambda_{r+1}/\hat\lambda_r = O_P(p^2n^{-1}/p^{2-2\delta}) =
O_P(p^{2\delta}n^{-1}) = O_P(h_n^2).
\]
This completes the proof of the corollary.
\end{pf*}

In the following, we use $\sigma_j(\bM)$ to denote the $j$th
largest singular value of a matrix $\bM$, so that $\sigma_1(\bM) =
\norm{\bM}$. We use $\lambda_j(\bM)$ to denote the $j$th largest
eigenvalue of $\bM$.
\begin{pf*}{Proof of Theorem~\ref{thm2}}
The first part of the theorem is actually Theorem 2 of~\cite{LYB11}.
We prove the other parts of the theorem.
From equation (22) of~\cite{LYB11},
the sample lag-$k$ autocovariance matrix for $\bve_t$ satisfies
%
\begin{equation}\label{app3}
\norm{\wh\bSigma_\epsilon(k)} = O_P(pn^{-1}).
\end{equation}

%

Note that~(\ref{app1}) together
with~(\ref{app3}) implies
\begin{eqnarray*}
\norm{\wh\bM- \bM} &=& O_P\bigl(p^{1-\delta}(p^{1-\delta/2}n^{-1/2} +
pn^{-1})\bigr) \\
&=& O_P\bigl(p^{2-2\delta}(p^{\delta/2}n^{-1/2} + p^{\delta
}n^{-1})\bigr)=O_P(p^{2-2\delta}\ell_n),
\end{eqnarray*}
since $\ell_n = p^{\delta/2}n^{-1/2} = o(1)$.
We also have $\norm{\hat\bB- \bB} = O_P(\ell_n)$, similar to the
proof of Theorem~\ref{thm1}.

With these, for $j=1,\ldots,r$, using decomposition~(\ref{appd1}), we have
\[
|\hat\lambda_j - \lambda_j| = O_P(\norm{\wh\bM- \bM}) =
O_P(p^{2-2\delta}\ell_n)
= O_P(p^{2-3\delta/2}n^{-1/2}),
\]
which is Theorem~\ref{thm2}(i).
For $j=r+1,\ldots,(k_0+1)r$, using decomposition~(\ref{appd2}), we
have
\begin{eqnarray*}
K_1 &=& O_P\bigl((p^{1-\delta/2}n^{-1/2} + pn^{-1})^2\bigr) = O_P(p^{2-\delta
}n^{-1} +
p^2n^{-2}) = O_P(p^{2-\delta}n^{-1}),\\
|K_2| &=& O_P(\norm{\wh\bM- \bM} \cdot\norm{\wh\bB- \bB}) =
O_P(p^{2-2\delta}\ell_n^2) = O_P(p^{2-\delta}n^{-1}),\\
|K_3| &=& O_P(\norm{\wh\bB- \bB}^2 \cdot\norm{\bM}) =
O_P(p^{2-2\delta}\ell_n^2) = O_P(p^{2-\delta}n^{-1}).
\end{eqnarray*}
Hence $\hat\lambda_j = O_P(p^{2-2\delta}\ell_n^2) =
O_P(p^{2-\delta}n^{-1})$, which is Theorem~\ref{thm2}(ii).

For part (iii), we define
\[
\bW_y(k_0) = (\bSigma_y(1), \ldots, \bSigma_y(k_0)),\qquad
\hat\bW_y(k_0) = (\hat\bSigma_y(1), \ldots, \hat\bSigma_y(k_0)),
\]
so that $\bM= \bW_y(k_0)\bW_y(k_0)'$ and $\hat\bM=
\hat\bW_y(k_0)\hat\bW_y(k_0)'$. We define similarly $\hat\bW_x(k_0),
\hat\bW_{x\epsilon}(k_0)$, $\hat\bW_{\epsilon x}(k_0)$ and
$\hat\bW_{\epsilon}(k_0)$. Then we can write
\[
\hat\bW_y(k_0) = M_1 + M_2 + \hat\bW_\epsilon(k_0),
\]
where $M_1 = \bA(\hat\bW_x(k_0)(\bI_{k_0} \otimes\bA') +
\hat\bW_{x\epsilon}(k_0))$, $M_2 = \hat\bW_{\epsilon
x}(k_0)(\bI_{k_0} \otimes\bA')$. It is easy to see that
\[
\operatorname{rank}(M_1) \leq r, \qquad \operatorname{rank}(M_2) \leq k_0r,
\]
so that rank$(M_1 + M_2) \leq(k_0+1)r$. This implies that
\[
\sigma_j(M_1+M_2) = 0\qquad \mbox{for } j=(k_0+1)r+1,\ldots,p.
\]
Then by Theorem 3.3.16(a) of~\cite{HJ91},
for
$j=(k_0+1)r + 1,\ldots,p$,
\begin{eqnarray*}
\hat\lambda_j &=& \lambda_j(\hat\bM) = \sigma_j^2(\hat\bW_y(k_0))
\leq \bigl(\sigma_j(M_1 + M_2) + \sigma_1(\hat\bW_\epsilon(k_0))\bigr)^2\\
&=& \sigma_1^2(\hat\bW_\epsilon(k_0)) \leq\sum_{k=1}^{k_0} \norm
{\wh\bSigma_\epsilon(k)}^2 = O_P(p^2n^{-2}),
\end{eqnarray*}
where the last equality sign follows from~(\ref{app3}). This proves
Theorem~\ref{thm2}(iii).

We prove Theorem~\ref{thm2}(ii)$'$ now. Using Lemma 3 of~\cite{LYB11},
with the same technique as in the proof of Theorem 1 in their paper, we
can write
%
\begin{equation}\label{app15}
\wh\bB= (\bB+ \bA\bP)(\bI+ \bP'\bP)^{-1/2}\qquad\mbox{with }
\norm{\bP} = O_P(\ell_n).
\end{equation}

With the definition of $\wh\bB$ as in the proof of Theorem~\ref{thm1}, we can write $\hat\lambda_{r+1}$,
the $(r+1)$th largest eigenvalue of $\wh\bM$, as the $(1,1)$ element
of\vspace*{1pt} the diagonal matrix $\wh\bD= \wh\bB'\wh\bM\wh\bB$, where
$\wh\bM\wh\bB= \wh\bB\wh\bD$. But from~(\ref{app15}), we also
have $\bB'\wh\bB= \bB'(\bB+ \bA\bP)(\bI+ \bP'\bP)^{-1/2} =
(\bI+ \bP'\bP)^{-1/2}$, hence
\[
(\bI+ \bP'\bP)^{1/2}\bB'\wh\bM\wh\bB= (\bI+ \bP'\bP
)^{1/2}\bB'\wh\bB\wh\bD= (\bI+ \bP'\bP)^{1/2}(\bI+ \bP'\bP
)^{-1/2}\wh\bD= \wh\bD.
\]
Further,\vspace*{1pt} by using Neumann series expansions of $(\bI+ \bP'\bP
)^{1/2}$ and $(\bI+ \bP'\bP)^{-1/2}$, we see that the largest order
term of $(\bI+ \bP'\bP)^{1/2}\bB'\wh\bM\wh\bB$ is contributed\vspace*{1pt}
from $\bB'\wh\bM(\bB+ \bA\bP)$, since from~(\ref{app15}) we have
$\norm{\bP} = O_P(\ell_n) = o_P(1)$. Hence\vspace*{1pt} the rate of $\hat\lambda
_{r+1}$ can be analyzed using the $(1,1)$ element of $\bB'\wh\bM(\bB
+ \bA\bP)$.

Some notation first. Define $\mathbf{1}_k$ the column vector of $k$
ones, and
\[
\bE_{r,s} = (\bve_r,\ldots,\bve_s), \qquad\bX_{r,s} = ({\mathbf
x}_r,\ldots,{\mathbf x}_s)\qquad\mbox{for } r \leq s.
\]
Since $k$ is finite and $\{\bve_t\}$ and $\{{\mathbf x}_t\}$ are
stationary, for convenience in this proof we take the sample lag-$k$
autocovariance matrix for $\{\bve_t\}$, $\{{\mathbf x}_t\}$ and the
cross lag-$k$ autocovariance matrix between $\{\bve_t\}$ and $\{
{\mathbf x}_t\}$ to be respectively, for $k > 0$,
\begin{eqnarray*}
\wh\bSigma_\epsilon(k) &=& n^{-1}\bigl(\bE_{k+1,n} -
(n-k)^{-1}\bE_{k+1,n}\mathbf{1}_{n-k}\mathbf{1}_{n-k}'\bigr)\\
&&{}\times\bigl(\bE_{1,n-k}
- (n-k)^{-1}\bE_{1,n-k}\mathbf{1}_{n-k}\mathbf{1}_{n-k}'\bigr)'\\
&=& n^{-1} \bE_{k+1,n}\mathbf{T}_{n-k}\bE_{1,n-k}',\\
\wh\bSigma_x(k) &=& n^{-1} \bX_{k+1,n}\mathbf{T}_{n-k}\bX_{1,n-k}'
\end{eqnarray*}
and
\[
\wh\bSigma_{x\epsilon}(k) = n^{-1} \bX
_{k+1,n}\mathbf{T}_{n-k}\bE_{1,n-k}',
\]
where $\mathbf{T}_j = \bI_j - j^{-1}\mathbf{1}_j\mathbf{1}_j'$. Then
\[
\bB'\wh\bM(\bB+ \bA\bP) = \sum_{k=1}^{k_0} \bB'\wh\bSigma
_y(k)\wh\bSigma_y(k)'(\bB+ \bA\bP) = \sum_{k=1}^{k_0} \bF_k(\bF
_k' + \bG_k),
\]
where
\begin{eqnarray*}
\bF_k &=&   n^{-1}\bB'\bE_{k+1,n}\mathbf
{T}_{n-k}\bX_{1,n-k}'\bA' + n^{-1}\bB'\bE_{k+1,n}\mathbf
{T}_{n-k}\bE_{1,n-k}',\\
\bG_k &=&   n^{-1}\bA\bX_{1,n-k}\mathbf{T}_{n-k}\bX_{k+1,n}'\bP'
+ n^{-1}\bE_{1,n-k}\mathbf{T}_{n-k}\bX_{k+1,n}'\bP'\\
&&{} + n^{-1}\bA\bX_{1,n-k}\mathbf{T}_{n-k}\bE_{k+1,n}'\bA\bP' +
n^{-1}\bE_{1,n-k}\mathbf{T}_{n-k}\bE_{k+1,n}'\bA\bP'.
\end{eqnarray*}
Some tedious algebra (omitted here) shows that the dominating term of
the above product is $\sum_{k=1}^{k_0}n^{-2}\bB'\bE_{k+1,n}\mathbf
{T}_{n-k}\bX_{1,n-k}'\bX_{1,n-k}\mathbf{T}_{n-k}\bX_{k+1,n}'\bP'$.
Defining ${\mathbf c}_{1,k}' = ({\mathbf a}_{r+1}'\bve_{k+1},\ldots
,{\mathbf a}_{r+1}'\bve_n)$ and $\mathbf{p}_1$ the first column of
$\bP'$, the $(1,1)$ element of the said term is then
\begin{eqnarray*}
&&\sum_{k=1}^{k_0}n^{-2}{\mathbf c}_{1,k}'\mathbf{T}_{n-k}\bX
_{1,n-k}'\bX_{1,n-k}\mathbf{T}_{n-k}\bX_{k+1,n}'\mathbf{p}_1\\
&&\qquad\leq\sum_{k=1}^{k_0}n^{-2} \norm{{\mathbf c}_{1,k}'}\norm{\mathbf
{p}_1}\norm{\mathbf{T}_{n-k}}^2\norm{\bX_{1,n-k}}^2 \norm{\bX
_{k+1,n}}\\
&&\qquad\leq4\sum_{k=1}^{k_0}\norm{n^{-1/2}{\mathbf c}_{1,k}}\norm{\bP
}\norm{n^{-1/2}\bX_{1,n-k}}^2\norm{n^{-1/2}\bX_{k+1,n}}\\
&&\qquad=O_P\bigl(\norm{n^{-1/2}{\mathbf c}_{1,1}} \cdot\ell_n \cdot
p^{(3-3\delta)/2}\bigr).
\end{eqnarray*}
In the last line we used $\norm{n^{-1/2}\bX_{1,n-k}} =
O_P(p^{(1-\delta)/2})$, by noting that
\begin{eqnarray*}
\norm{n^{-1/2}\bX_{1,n-k}}^2 &=& \norm{n^{-1}\bX_{1,n-k}\bX
_{1,n-k}'} \asymp\norm{n^{-1}\bX_{1,n-k}\mathbf{T}_{n-k}\bX
_{1,n-k}'}\\
&=& \norm{\hat\bSigma_x(0)} \leq\norm{\hat\bSigma_x(0) - \bSigma
_x(0)} + \norm{\bSigma_x(0)}\\
&=& O_P(p^{1-\delta}n{-1/2}) + O_P(p^{1-\delta}) = O(p^{1-\delta}),
\end{eqnarray*}
where $\norm{\hat\bSigma_x(0) - \bSigma_x(0)} = O_P(p^{1-\delta
}n{-1/2})$ is from~(\ref{app16}) and $\norm{\bSigma_x(0)} =
O(p^{1-\delta})$ is assumed in condition (C5).
With condition (C9), we can show that $\norm{n^{-1/2}{\mathbf c}_{1,1}}
= O_P(1)$, since
\begin{eqnarray*}
P(\norm{n^{-1/2}{\mathbf c}_{1,1}} > x) &=& P\Biggl(n^{-1}\sum
_{j=k+1}^n{\mathbf a}_{r+1}'\bve_j\bve_j'{\mathbf a}_{r+1} >
x^2\Biggr)\\
&\leq& (n-k){\mathbf a}_{r+1}'\bSigma_{\bve}{\mathbf a}_{r+1}/(nx^2)
\leq\sigma_{\max}^2 / x^2,
\end{eqnarray*}
where we used the Markov inequality with $\sigma_{\max}^2$ the
maximum eigenvalue of~$\bSigma_{\bve}$,
and the fact that ${\mathbf a}_{r+1}'{\mathbf a}_{r+1} = 1$. Hence the
$(1,1)$ element of $\bB'\wh\bM(\bB+ \bA\bP)$
has rate $O_P(p^{(3-3\delta)/2}\ell_n) = O_P(p^{3/2 - \delta
}n^{-1/2})$, which is also the rate of~$\hat\lambda_{j}$ for $j \geq r+1$.
This completes the
proof of the theorem.
\end{pf*}

We outline the proofs of Theorems~\ref{thm3} and~\ref{thm4}
below. Detailed proofs can be found in the supplemental article (Lam
and Yao~\cite{LY12}).
%
\begin{pf*}{Outline proof of Theorem~\ref{thm3}}
First, under model (\ref {e6}) and $\bM$ defined in~(\ref{d3}), with
conditions (C1)--(C4), (C5)$'$, (C6)$'$, we can show that the rates of the
eigenvalues of $\bM$ are given by
%
\begin{equation}\label{appo1}\qquad
\lambda_j \asymp
\cases{
p^2, &\quad for $j=1,\ldots,r_1$;\vspace*{2pt}\cr
p^{2-2\delta_2}, &\quad if $\norm{\bW_{21}}_{\min} = o(p^{1-\delta_2})$
$(c=1)$\vspace*{2pt}\cr
&\quad for $j=r_1+1,\ldots,r$;\vspace*{2pt}\cr
p^{2-2c\delta_2}, &\quad if $\norm{\bW_{21}}_{\min}
\asymp p^{1-c\delta_2}$, $1/2 \leq c < 1$, and\vspace*{2pt}\cr
&\quad $\norm{\bW_1\bW_{21}'} \leq q \norm{\bW_1}_{\min} \norm{\bW_{21}},
0 \leq q < 1$,\vspace*{2pt}\cr
&\quad for $j=r_1+1,\ldots,r$.}
\end{equation}
For model~(\ref{e6}), and $\bM^*$ defined in Section~\ref{sec6} by
$\by_t^*$ in~(\ref{e7}), we have
%
\begin{equation}\label{appo2}
\lambda_j^* \asymp p^{2-2\delta_2}
\qquad\mbox{for } j=1,\ldots,r_2.
\end{equation}

We cannot use Lemma 3 of~\cite{LYB11}
to prove this theorem for the one-step estimation, since the condition
$\norm{\bE} \leq\operatorname{sep}(\bD_1,\bD_2)/5$ gives a restrictive
condition on the growth rate of $p$, and also restricts the range of
$\delta_2$ allowed. Instead, we use Theorem 4.1 of~\cite{S73}.

Write $\bM= \bX_{ij}\bD_{ij}\bX_{ij}'$ for $i \neq j =1,2$, where
$\bX_{ij} = (\bA_i   \bA_j   \bB)$, $\bB$ is the orthogonal
complement of $\bA= (\bA_1   \bA_2)$, and $\bD_{ij}$ is diagonal
with $\bD_{ij} = \operatorname{diag}(\bD_i,   \bD_j,   \bzero)$
where $\bD_1$ contains $\lambda_j$ for $j=1,\ldots,r_1$ and $\bD_2$
contains $\lambda_j$ for $j=r_1+1,\ldots,r$. With $\bE= \wh\bM-
\bM$, define
\[
\bX'\bE\bX= (\bE_{ij})   \qquad\mbox{for } {1 \leq i,j \leq3},
\]
where $\bE_{ij} = \bA_i'\bE\bA_j$ if we denote $\bB= \bA_3$.

Define $\operatorname{sep}(\bM_1, \bM_2) = {\min_{\lambda\in\lambda(\bM
_1), \mu\in\lambda(\bM_2)}}|\lambda- \mu|$. If we can show that
%
\begin{eqnarray}\label{appo3}
\norm{(\bE_{ij}, \bE_{i3})} = o_P(\gamma_{ij})\qquad\qquad \nonumber\\[-8pt]\\[-8pt]
&&\eqntext{\mbox{with }
\gamma_{ij} = \operatorname{sep}\left( \bD_i + \bE_{ii},
\pmatrix{
\bD_j + \bE_{jj} & \bE_{j3} \cr
\bE_{3j} & \bE_{33}}
\right),}
\end{eqnarray}
then condition (4.2) in~\cite{S73}
is satisfied asymptotically, so that we can use their Theorem 4.1 to
conclude that for $i\neq j = 1,2$,
%
\begin{equation}\label{appo4}
\norm{\wh\bA_i - \bA_i} = O_P\bigl(\norm{(\bE_{ij}, \bE_{i3})}/\gamma_{ij}\bigr).
\end{equation}
Since we can show that $\norm{\bE_{12}} = O_P(\norm{\bE_{13}}) =
O_P(p^2n^{-1/2})$, we have $\norm{(\bE_{12}$, $\bE_{13})} = O_P(p^2n^{-1/2})$.
We can also show that $\gamma_{12} \asymp p^2$ using~(\ref{appo1}).
Hence (\ref{appo3}) is satisfied, and~(\ref{appo4}) implies that
\[
\norm{\wh\bA_1 - \bA_1} = O_P(p^2n^{-1/2} / p^2) = O_P(n^{-1/2}).
\]
Also, we can show that $\norm{\bE_{23}} = O_P(\bE_{21}) =
O_P(p^{2-\delta_2/2}n^{-1/2})$, implying that $\norm{(\bE_{21}, \bE
_{23})} = O_P(p^{2-\delta_2/2}n^{-1/2})$. We can also show that
$\gamma_{21} \asymp p^{2-2c\delta_2}$ using~(\ref{appo1}), provided
$p^{c\delta_2}n^{-1/2} \to0$. Hence~(\ref{appo3}) is satisfied since
we assumed $\nu_n \to0$, and so~(\ref{appo4}) implies that
\[
\norm{\wh\bA_2 - \bA_2} = O_P(p^{2-\delta_2/2}n^{-1/2} /
p^{2-2c\delta_2}) = O_P\bigl(p^{(2c-1/2)\delta_2}n^{-1/2}\bigr) = O_P(\nu_n),
\]
which completes the proof for the one-step estimation.

For the two-step estimation, write $\bM^* = (\bA_2   \bB^*)\bD
^*(\bA_2   \bB^*)'$, where $\bB^*$ is the orthogonal complement of
$\bA_2$, and $\bD^*$ is diagonal with $\bD^* = \operatorname{diag}(\bD
_2^*,\mathbf{0})$. The matrix $\bD_2^*$ contains $\lambda_j^*$ for
$j=1,\ldots,r_2$, so that~(\ref{appo2}) implies $\operatorname{sep}(\bD
_2^*, \mathbf{0}) \asymp p^{2-2\delta_2}$.

We can show that $\norm{\bE^*} = \norm{\wh\bM^* - \bM^*} =
O_P(p^{2-2\delta_2}\kappa_n)$, hence $\norm{\bE^*} = o_P(
\operatorname{sep}(\bD_2^*,\mathbf{0}))$, as $\kappa_n \to0$. Hence we can use
Lemma 3 of~\cite{LYB11}
to conclude that
\[
\norm{\wt\bA_2 - \bA_2} = O_P\bigl(\norm{\bE_{21}^*}/\operatorname{sep}(\bD
_2^*, \mathbf{0})\bigr).
\]
Since we can show that $\norm{\bE_{21}^*} = O_P(p^{2-3\delta
_2/2}n^{-1/2})$, we then have
\[
\norm{\wt\bA_2 - \bA_2} = O_P(p^{2-3\delta
_2/2}n^{-1/2}/p^{2-2\delta_2}) = O_P(p^{\delta_2/2}n^{-1/2}),
\]
which completes the proof of the theorem.
\end{pf*}

To prove Theorem~\ref{thm4}, we need two lemmas first.
%
\begin{lemma}\label{lemma}
Under the same conditions and notations of Theorem~\ref{thm3}, the
following assertions hold:

\begin{longlist}[(iii)$'$]
\item[(i)] For $j=1,\ldots,r_1$, $|\hat\lambda_j - \lambda_j| =
O_P(p^{2}n^{-1/2}) $.

\item[(ii)] For $ j=r_1+1,\ldots,r$, $|\hat\lambda_j - \lambda_j| =
O_P(p^{2}(n^{-1/2} + \nu_{n}^2)) $
provided $\nu_n \to0$, $p^{c\delta_2}n^{-1/2} \to0$.

\item[(iii)] For $j=r+1,\ldots,p$, $\hat\lambda_j = O_P(p^{2}\nu_{n}^2)$,
provided $\nu_n \to0$, $p^{c\delta_2}n^{-1/2} \to0$.

\item[(iv)] For $j=1,\ldots,r_2$, $|\hat\lambda_j^* - \lambda_j^*| =
O_P(p^{2-2\delta_2}\kappa_{n})$.

\item[(v)] For $j=r_2+1,\ldots,p$, $\hat\lambda_j^* = O_P(p^{2-2\delta
_2}\kappa_{n}^2)$.

\item[(vi)] For $j=(k_0+1)r + 1,\ldots,p$, $ \hat\lambda_j,   \hat\lambda
_j^* = O_P(p^2n^{-2}) =
O_P(p^{2-2\delta_2}\kappa_{n}^4)$.

\item[(iii)$'$] If in addition condition \textup{(C9)} holds, then for $j=r+1,\ldots,p$,
$\hat\lambda_j = O_P(p^{3/2}\nu_{n})$, provided $\nu_n \to0$,
$p^{c\delta_2}n^{-1/2} \to0$.
\end{longlist}
\end{lemma}

The proof of this lemma is left in the supplementary materials for this
paper. Together with~(\ref{appo1}) and~(\ref{appo2}), we have the
following lemma.\vadjust{\goodbreak}
%
\begin{lemma}
\label{lemma2}
Let conditions \textup{(C1)--(C4)}, \textup{(C5)$'$},
\textup{(C6)$'$}, \textup{(C7)} and \textup{(C8)} hold. Then as $n,p \to
\infty$ with $n = O(p)$, and with $\nu_{n}$, $\kappa_{n} \to0$ the
same as in Theorem~\ref{thm3} and $p^{c\delta_2}n^{-1/2} \to0$, we have
\[
\hat\lambda_{j+1}/\hat\lambda_j \cases{
\asymp1, &\quad $j=1,\ldots,r_1-1$;\vspace*{2pt}\cr
= O_P(n^{-1/2} + \nu_{n}^2 + p^{- 2\delta_2}),
&\quad $j=r_1 $, if $\norm{\bW_{21}}_{\min} = o(p^{1-\delta_2})$\vspace*{2pt}\cr
&\quad $(c=1)$;\vspace*{2pt}\cr
= O_P(n^{-1/2} + \nu_{n}^2 + p^{- 2c\delta_2}),
&\quad $j=r_1$, if $\norm{\bW_{21}}_{\min} \asymp p^{1-c\delta_2}$\vspace*{2pt}\cr
&\quad for $1/2 \leq c < 1$, and\vspace*{2pt}\cr
&\quad$\norm{\bW_1\bW_{21}'} \leq q \norm{\bW_1}_{\min} \norm{\bW
_{21}} $\vspace*{2pt}\cr
&\quad for $0 \leq q <1 $.}
\]
Furthermore, if $\norm{\bW_{21}}_{\min} = o(p^{1-\delta_2})$ and
$p^{5\delta_2/2}n^{-1/2} \to0$, we have
\[
\hat\lambda_{j+1}/\hat\lambda_j \cases{
\asymp1, &\quad $j=r_1 + 1,\ldots,r-1$;\vspace*{2pt}\cr
= O_P(p^{2\delta_2}\nu_{n}^2), &\quad $j=r$;\vspace*{2pt}\cr
= O_P(p^{2\delta_2 - 1/2}\nu_n), &\quad $j=r$, and condition \textup{(C9)} holds.}
\]

If $\norm{\bW_{21}}_{\min} \asymp p^{1-c\delta_2}$ for $1/2 \leq c <
1$, $ \norm{\bW_1\bW_{21}'} \leq q \norm{\bW_1}_{\min} \norm{\bW_{21}}
$ for $0 \leq q <1 $, and $p^{(3c-1/2)\delta_2}n^{-1/2} \to0$, we have
\[
\hat\lambda_{j+1}/\hat\lambda_j \cases{
\asymp1, &\quad $j=r_1 + 1,\ldots,r-1$;\cr
= O_P(p^{2c\delta_2}\nu_{n}^2), &\quad $ j=r $;\vspace*{1pt}\cr
= O_P(p^{2c\delta_2 - 1/2}\nu_n), &\quad $j=r$, and condition \textup{(C9)} holds.}
\]

For the two-step procedure, let conditions \textup{(C1)--(C4)},
\textup{(C5)$'$}, \textup{(C6)$'$}, \textup{(C7)} and \textup{(C8)}
hold and $n = O(p)$. Then we have
\[
\hat\lambda_{j+1}^*/\hat\lambda_j^* \cases{
\asymp1, &\quad $j=1,\ldots,r_2-1$;\cr
= O_P(\kappa_{n}^2), &\quad $ j=r_2 $.}
\]
\end{lemma}
\begin{pf}
We only need to find the asymptotic
rate for each $\hat\lambda_j$ and $\hat\lambda_j^*$. The rate of
each ratio can then be obtained from the results of Lemma~\ref{lemma}.

For\vspace*{1pt} $j=1,\ldots,r_1$, from Lemma~\ref{lemma}, $\norm{\hat\lambda_j
- \lambda_j} = O_P(p^2n^{-1/2}) = o_P(\lambda_j)$, and hence $\hat
\lambda_j \asymp\lambda_j \asymp p^2$ from~(\ref{appo1}).

Consider\vspace*{1pt} the case $\norm{\bW_{21}}_{\min} \asymp p^{1-c\delta_2}$.
For $j=r_1+1,\ldots,r$, since $|\hat\lambda_j - \lambda_j| =
O_P(p^2(n^{-1/2} + \nu_n^2))$, we have $\hat\lambda_j \leq\lambda
_j + O_P(p^2(n^{-1/2} + \nu_n^2)) = O_P(p^{2-2c\delta_2} + p^2\nu
_n^2 + p^2n^{-1/2})$, and hence
\[
\hat\lambda_{r_1+1}/\hat\lambda_{r_1} = O_P\bigl((p^{2-2c\delta_2} +
p^2\nu_n^2 + p^2n^{-1/2}) / p^2\bigr) = O_P(n^{-1/2} + \nu_n^2 +
p^{-2c\delta_2}).
\]
The other case is proved similarly.

For $j=r_1+1,\ldots,r$, to make sure $\hat\lambda_j$ will not be
zero or close to zero, we need
\[
|\hat\lambda_j - \lambda_j| = O_P\bigl(p^2(n^{-1/2} + \nu_n^2)\bigr) =
o_P(\lambda_j),
\]
where $\lambda_j \asymp p^{2-2c\delta_2}$ as in~(\ref{appo1}). Hence
we need $p^2(n^{-1/2} + \nu_n^2) = o(p^{2-2c\delta_2})$, which is
equivalent to the condition $p^{(3c-1/2)\delta_2}n^{-1/2} \to0$. With
this condition satisfied, then $\hat\lambda_j \asymp\lambda_j
\asymp p^{2-2c\delta_2}$ for $j=r_1+1,\ldots,r$. Since $\hat\lambda
_j = O_P(p^2\nu_n^2)$ for $j=r+1,\ldots,p$, we then have
\[
\hat\lambda_{r+1}/\hat\lambda_r = O_P(p^2\nu_n^2 / p^{2-2c\delta
_2}) = O_P(p^{2c\delta_2}\nu_n^2).
\]
All other rates can be proved similarly, and thus are omitted.
\end{pf}
\begin{pf*}{Proof of Theorem~\ref{thm4}}
With Lemma~\ref{lemma2}, Theorem~\ref{thm4}(i) is obvious. For
Theorem~\ref{thm4}(ii), note that the range of $\delta_2$ and the
rates given in the theorem ensure that $n^{-1/2} + \nu_n^2 +
p^{-2c\delta_2} = o(p^{2c\delta_2-1/2}\nu_n) = o(p^{2c\delta_2}\nu
_n^2)$. Hence Lemma~\ref{lemma2} implies a better rate of convergence
for $\hat\lambda_{r_1+1}/\hat\lambda_{r_1}$ no matter whether
condition (C9) holds or not.
We can use a similar argument to prove part (iii), and details are thus
omitted.
\end{pf*}

\end{appendix}

\section*{Acknowledgments}

We are grateful to the Joint Editor Professor Peter B\"uhlmann, the
Associate Editor and the two referees for their helpful comments and
suggestions.

\begin{supplement}
\stitle{Detailed proofs of Theorems~\ref{thm3} and~\ref{thm4}}
\slink[doi]{10.1214/12-AOS970SUPP} 
\sdatatype{.pdf}
\sfilename{aos970\_supp.pdf}
\sdescription{The document contains detailed proofs of Theorem~\ref{thm3} and
\ref{thm4} in the paper.}
\end{supplement}


\printaddresses

\end{document}